\documentclass{amsart}

\usepackage{amssymb,latexsym,amsfonts,amsmath}
\usepackage{graphicx,psfrag,epsfig,epsf}
\include{diagrams}
\usepackage{subfigure}
\usepackage{graphicx}
\usepackage{color}
\usepackage{diagrams}
\usepackage{algorithm2e}
\usepackage{tikz}
\usetikzlibrary{automata}   

\usepackage{amsfonts}

\newcommand{\alt}{\mathrm{alt}}

\usepackage{amssymb,latexsym,amsfonts,amsmath}
\usepackage{graphicx}

\newtheorem{theorem}{Theorem}[section]
\newtheorem{lemma}[theorem]{Lemma}

\newtheorem{proposition}[theorem]{Proposition}

\theoremstyle{definition}
\newtheorem{definition}[theorem]{Definition}

\theoremstyle{remark}
\newtheorem{remark}[theorem]{Remark}
\numberwithin{equation}{section}

\begin{document}

\begin{abstract}
Time--delay systems are an important class of dynamical systems that provide a solid mathematical framework to deal with many application domains of interest. In this paper we focus on nonlinear control systems with unknown and time--varying delay signals and we propose one approach to the control design of such systems, which is based on the construction of symbolic models. Symbolic models are abstract descriptions of dynamical systems where one symbolic state and one symbolic input correspond to an aggregate of states and an aggregate of inputs.
We first introduce the notion of incremental input--delay--to--state stability and characterize it by means of Lyapunov--Krasovskii functionals.  We then derive sufficient conditions for the existence of symbolic models that are shown to be alternating approximately bisimilar to the original system. Further results are also derived which prove the computability of the proposed symbolic models in a finite number of steps.
\end{abstract}

\title[Symbolic Models for Time--Varying Time--Delay Systems]{Symbolic Models for Nonlinear Time--Varying Time--Delay Systems\\ via Alternating Approximate Bisimulation}
\thanks{This work has been partially supported by the Center of Excellence for Research DEWS and by the Network of Excellence HyCon$^{2}$, Grant agreement No. 257462.}

\author[Giordano Pola, Pierdomenico Pepe and Maria D. Di Benedetto]{
Giordano Pola$^{1}$, Pierdomenico Pepe$^{1}$ and Maria D. Di Benedetto$^{1}$}
\address{$^{1}$
Department of Electrical and Information Engineering, Center of Excellence DEWS,
University of L{'}Aquila, Poggio di Roio, 67040 L{'}Aquila, Italy}
\email{ \{giordano.pola,pierdomenico.pepe,mariadomenica.dibenedetto\}@univaq.it}
\urladdr{
http://www.diel.univaq.it/people/pola/
}
\urladdr{
http://www.diel.univaq.it/people/pepe/
}
\urladdr{
http://www.diel.univaq.it/people/dibenedetto/
}
%

\maketitle

\vspace{5mm}
\textbf{keywords: } 
Time--delay systems, symbolic models, alternating approximate bisimulation, incremental input--delay--to--state stability, time--varying delays.

\section{Introduction}

Time--delay systems are an important class of dynamical systems that provide a sound mathematical framework to deal with many application domains of interest, from biology, to chemical, electrical, and mechanical engineering, and economics (see e.g. \cite{Niculescu:01}). 
Over the years, several researchers focused on stability, stabilization, regulation and linearization problems for time--delay systems (see e.g. \cite{Hale:93,Niculescu:01,Kolmanovskii:99, surveyTDS} and the references therein). 
However, the increasing complexity of current technology asks for control algorithms that can 
also deal with a different kind of specifications, such as 
safety and liveness properties, obstacle avoidance, fairness constraints, language and logic specifications, among many others (see e.g. \cite{LTLControl,paulo}). 
In this paper, we consider nonlinear control systems with unknown and time--varying delay signals and we propose a method to deal with such specifications which is based on symbolic models. 
Symbolic models are abstract models where each symbolic state and each symbolic label represent an aggregation of continuous states and an aggregation of input and delay signals in the original model. 
Since these symbolic models are of the same nature as the models used in computer science to describe
software and hardware, they provide a unified language to study problems of control where software and hardware interact with the physical world.
Moreover, the use of symbolic models allows one to leverage the rich literature developed in the computer science community,
as for example supervisory control \cite{RamWonI} and algorithmic game theory \cite{AVW03}, for control design of purely continuous processes.
Many researchers have recently faced the problem of identifying classes of control systems admitting symbolic models.
For example, controllable linear control systems and incrementally stable nonlinear control systems were shown in \cite{LTLControl} and respectively in \cite{PolaAutom2008,PolaSIAM2009}, to admit symbolic models. Symbolic models for multi--affine systems have been proposed in \cite{Belta:06} and for nonlinear switched systems in \cite{GirardTAC2010}. 
Symbolic models for nonlinear systems with known and constant delays have been explored in \cite{PolaSCL10}. In this paper we extend the results of the work in \cite{PolaSCL10} to nonlinear control systems with unknown and time--varying delays.
This class of time--delay systems arises in modeling many application domains of interest and are also challenging from a theoretical point of view, because of their inherent complexity, see e.g. \cite{Fridman06} and the references therein. 
We stress that the extension of a control theory for systems described by time--delay functional differential equations, from the case of known time--delays to the case of unknown time--delays, is in general not straightforward, 
as well as the extension from the case of constant time--delays to the case of time--varying time--delays. As far as the stabilization of general systems with time--varying time--delays is concerned, by the use of control Lyapunov functionals, the reader can refer to the recent book \cite{KarJiaBook}. The methodology of control Lyapunov functionals, in general, requires the knowledge of the time--delays. Most nonlinear control laws available in the literature, mainly dealing with the stabilization and the input--output linearization problems, are in general designed under the assumption of knowledge of the time--delay (often assumed constant), see e.g. \cite{Jankovic:01,GermaniSIAM00,Oguchi:02,IOMM04}. A big effort is being made by researchers in recent years to deal with unknown or time--varying time--delays, for instance to manage networking control, where the communication time--delay is often time--varying and unknown, while an estimation of its upper and lower bounds can generally be given with a certain confidence (see for instance \cite{GeAut05} for an adaptive controller in the case of unknown constant delays). It is shown in \cite{KarJiaBook,KJESAIM10} that suitable triangular time--delay systems can be stabilized by a delay-free controller, thus completely solving the stabilization problem for this class of systems, also in the case of unknown time--varying time--delays. However, to our knowledge, the control of general nonlinear systems with unknown time--varying time--delays is still an open and challenging research problem in many cases, for instance when the above mentioned specifications are required. \\
The framework considered in this paper requires the results of \cite{PolaSCL10} to be extended in several directions. First, an appropriate generalization of the notion of incremental stability proposed in \cite{PolaSCL10} is needed, in order to quantify the mismatch of state trajectories in the symbolic models and in the \textit{time--varying} time--delay systems. To this aim the notion of \textit{incremental input--delay--to--state stability} ($\delta$--IDSS) is here introduced  and characterized in terms of Lyapunov--Krasovskii functionals. The \textit{incremental input--delay--to--state stability} ($\delta$--IDSS) is a novel notion that characterizes the relationship between solutions obtained by different time--delay signals, besides by different initial conditions and inputs. Second, the presence of the additional time--varying time--delay signal in the differential equation requires appropriate approximating schemes in order to incorporate its effects in the symbolic model. We solve this problem by approximating this time--delay signal by a (quantized) first--order spline \cite{SplineBook}, that then can be regarded as an additional unknown disturbance label in the symbolic model. 
To our knowledge, it is the first time that the (unknown) time--delay signal is approximated by piece-wise linear functions, for the purpose of building up a suitable controller. Generally, for this aim, the only state space is projected into suitable finite dimensional subspaces (see \cite{GermaniSIAM00} and references therein). 
Third, an appropriate notion of approximate equivalence is needed in order to capture the adversarial nature of the control input labels and the aforementioned disturbance labels. To this aim, we resort to the notion of alternating approximate bisimulation, recently introduced in \cite{PolaSIAM2009}, that guarantees that control strategies synthesized on the symbolic models, based on alternating approximate bisimulations, can be readily transferred to the original model, independently of the particular realization of the time--varying delay signals. While the work in \cite{PolaSCL10} shows the existence of symbolic models that are \textit{approximately bisimilar} \cite{AB-TAC07} to \textit{incrementally input--to--state stable} nonlinear control systems with \textit{constant} and \textit{known} time--delays, the results in the present work show existence of symbolic models that are \textit{alternating approximately bisimilar} \cite{PolaSIAM2009} to \textit{incrementally--input--delay--to--state stable} nonlinear control systems with \textit{time--varying} and \textit{unknown} time--delays. Further results are also derived which prove the computability of the proposed symbolic models in a finite number of steps.
An illustrative example is finally presented which shows the control design of a $2$--dimensional time--delay system with a synchronization--type specification. 
A preliminary investigation on the existence of symbolic models for time--varying time--delay systems appeared in the conference publication \cite{PolaCDC2010}. In this paper we present a detailed and mature description of the results announced in \cite{PolaCDC2010}, including proofs and an example. \\
The paper is organized as follows. In Section 2 we introduce the class of time--delay systems under consideration. In Section 3 we introduce and characterize the notion of incremental--input--delay--to--state stability. In Section 4 we recall the notions of approximate (bi--)simulations and of alternating approximate (bi--)simulations. Section 5 is devoted to the study of existence of symbolic models and Section 6 provides constructive results of the proposed symbolic models. In Section 7 we present an illustrative example. Section 8 offers some concluding remarks. 
For the sake of completeness, a detailed list of the employed notation is included in the Appendix (Section \ref{sec:notation}).

\section{Time--Varying Time--Delay Systems}

In this paper we consider the following nonlinear time--varying time--delay system:
\begin{equation}
\label{TDS}
\left\{
\begin{array}{ll}
\dot{x}(t)=f(x(t),x(t-\Delta(t)),u(t-r)),& t \in \mathbb{R}_0^{+}, a.e. \\
x(t)=\xi_{0}(t), & t \in [-\Delta_{\max},0],
\end{array}
\right.
\end{equation}
where:

\begin{itemize}

\item $x(t)\in \mathbb{R}^{n}$ and $x_{t}\in C^{0}([-\Delta_{\max},0];\mathbb{R}^{n})$ is the state at time $t\in\mathbb{R}_0^{+}$; we recall that $x_t(\theta)=x(t+\theta)$ with $\theta \in [-\Delta_{\max},0]$ and $t\in \mathbb{R}^{+}_{0}$;

\item $\xi_{0}=x_0\in C^{0}([-\Delta_{\max},0];\mathbb{R}^{n})$ is the initial condition;

\item $u(t)\in \mathbb{R}^{m}$ is the control input at time $t\in [-r,+\infty [$ and $r\in\mathbb{R}^{+}_{0}$ is the constant control input delay;

\item $\Delta:\mathbb{R}_{0}^{+}\rightarrow [\Delta_{\min},\Delta_{\max}]$ is the unknown time--varying delay with \mbox{$\Delta_{\min},\Delta_{\max}\in \mathbb{R}^{+}_{0}$};

\item $f: \mathbb {R}^n\times \mathbb {R}^n\times \mathbb {R}^m \rightarrow \mathbb{R}^{n}$ is the vector field.

\end{itemize}
\bigskip
\textbf{Assumption (A.0).}
In this paper we suppose that:
\begin{itemize}
\item Function $f$ satisfies $f(0,0,0)=0$ and it is Lipschitz on bounded sets, i.e. for every bounded set \mbox{$K\subset \mathbb {R}^n\times \mathbb {R}^n\times \mathbb {R}^m$}, there exists a constant $\kappa \in\mathbb{R}^{+}$ such that
\[
\begin{array}{c}
\Vert f(x_{1},y_{1},u_{1})  -f(x_{2},y_{2},u_{2})\Vert  \leq 
\kappa(\Vert x_{1}-x_{2}\Vert + \Vert
y_{1}-y_{2}\Vert + \Vert u_{1}-u_{2} \Vert),
\nonumber
\end{array}
\]
for all $(x_{1},y_1,u_{1}),(x_{2},y_2,u_{2})\in K$;

\item $\Vert u(t)\Vert \le B_{U}$, $\forall t \in [-r,+\infty[$ for some known $B_U \in \mathbb{R}^{+}$;

\item Delay bounds $r$, $\Delta_{\min}$ and $\Delta_{\max}$ are known;

\item The initial condition $\xi_0\in C^1([-\Delta_{\max},0];\mathbb{R}^n)$ with
\begin{eqnarray}
\begin{array} {cc}
\sup_{\tau \in [-\Delta_{\max},0]}\Vert\xi_0(\tau)\Vert\le B_{X}^{0}, &
\sup_{\tau \in[-\Delta_{\max},0]}\Vert\dot \xi_0(\tau)\Vert\le M_1,
\end{array}
\end{eqnarray}
where $B_{X}^{0}\in\mathbb{R}^{+}$ and $M_1 \in \mathbb{R}^{+}_{0}$ are known;

\item function $\Delta\in \mathcal{D}$ is unknown and $\mathcal{D}$ is the set of all continuously differentiable functions
$\Delta: \mathbb{R}^{+}_{0}\rightarrow
[\Delta_{\min},\Delta_{\max}]$ with $\Vert \dot \Delta(t) \Vert
\le d_{\min}$, where $d_{\min}\in[0,1[$ is known.
\end{itemize}
\bigskip
We denote by $U$ the set $\mathcal {B}_{B_U}(0)$ and by $\mathcal{U}$ the class of measurable control input signals from $[-r,+\infty[$ to $U$. Moreover the symbols $x(t,\xi_0,u,\Delta)$ and $x_t(\xi_0,u,\Delta)$ denote the solution at time $t\in\mathbb{R}^{+}_{0}$ in $\mathbb {R}^n$ and respectively in $C^0([-\Delta_{\max},0];\mathbb {R}^n)$ of the time--delay system in (\ref{TDS}) with initial condition $\xi_0$, input $u\in \mathcal{U}$ and time--delay signal $\Delta\in\mathcal{D}$.
In the sequel we refer to a time--delay system as in (\ref{TDS}) and satisfying Assumption (A.0), by means of the tuple:
\begin{equation}
\Sigma=(\mathbb{R}^{n},C^{0}([-\Delta_{\max},0];\mathbb{R}^{n}),\xi_{0},U,\mathcal{U},\mathcal{D},f),
\label{TupleTDS}
\end{equation}
where each entity has been defined before.

\section{Incremental Input--Delay--to--State Stability}\label{sec:stab}

The results of this paper rely upon some stability notions that
we introduce and characterize in this section. We start by recalling the notion of input--to--state stability.

\begin{definition}
\label{def:ISS}
\cite{SontagPrincipe,PepeJiangSCL}
A time--delay system $\Sigma$ is Input--to--State Stable (ISS) (uniformly with respect to the time--delay function) if there exist a $\mathcal{KL}$
function $\beta_{ISS}$ and a $\mathcal{K}$ function $\gamma_{ISS}$
such that for any time $t\in\mathbb{R}^{+}_{0}$, any initial
condition \mbox{$\xi_{0}\in
C^{0}([-\Delta_{\max},0];\mathbb{R}^{n})$}, any input $u\in
\mathcal{U}$, and any time--delay function $\Delta\in \mathcal{D}$, the solution $x_{t}(\xi_{0},u,\Delta)$ of (\ref{TDS}) exists for any $t\geq 0$ and, furthermore, the following inequality holds:
\begin{eqnarray}
\label{ineqiss}
\begin{array}{rcl}
\left\Vert x_{t}(\xi_{0},u,\Delta)\right\Vert_{\infty} & \leq & \max\{\beta_{ISS}(\left\Vert \xi_{0}\right\Vert_{\infty},t),\gamma_{ISS}(\left\Vert u|_{[-r,t-r[} \right\Vert
_{\infty})\}.
\end{array}
\end{eqnarray}
\end{definition}

A characterization of the ISS property in terms of Lyapunov--Krasovskii functionals can be found in \cite{PepeJiangSCL}. 

In the sequel, we suppose that:

\bigskip
\textbf{Assumption (A.1).} The time--delay system in (\ref{TupleTDS}) is ISS.
\bigskip

Under the above assumption any solution $x(t,\xi_0,u,\Delta)$ (the initial condition is here considered as
part of the solution) belongs to $\mathcal{B}_{B_{X}}(0)$ for any time $t\geq -\Delta_{\max}$, with $B_X=\max \{\beta_{ISS}(B_{X}^{0},0),\gamma_{ISS}(B_{U})\}$, where $\beta_{ISS}$ and $\gamma_{ISS}$ are the $\mathcal{KL}$ and respectively, the $\mathcal{K}$ function appearing in (\ref{ineqiss}). We denote by $X$ the set $\mathcal B_{B_X}(0)$ and by $\mathcal{X}$ the set $C^0([-\Delta_{\max},0];X)$. Under Assumption (A.1) there is no loss of generality in replacing $\mathbb{R}^{n}$ and $C^{0}([-\Delta_{\max},0];\mathbb{R}^{n})$ in (\ref{TupleTDS}), by $X$ and $\mathcal{X}$ respectively, resulting in:
\begin{equation}
\Sigma=(X,\mathcal{X},\xi_{0},U,\mathcal{U},\mathcal{D},f).
\label{TTDS}
\end{equation}

We can now introduce the notion of incremental input--delay--to--state stability, that adapts the notion of incremental input--to--state--stability, introduced in \cite{incrementalS} for the class of nonlinear control systems, to the class of time--varying time--delay systems considered in this paper.

\begin{definition}\label{definitiondeltaiss}
A time--varying delay system $\Sigma$ satisfying Assumptions (A.0) and (A.1), is \textit{incrementally Input--Delay--to--State Stable} (\mbox{$\delta$--IDSS}) if there exist a $\mathcal{KL}$ function $\beta$ and $\mathcal{K}$ functions $\gamma_{U}$, $\gamma_{D}$ such that for any time $t\in\mathbb{R}^{+}_{0}$, any initial conditions $\xi_{1},\xi_{2}$, any inputs $u_{1},u_{2}$ and any time--delay functions $\Delta_1,\Delta_2$, the corresponding solutions $x_{t}(\xi_{1},u_{1},\Delta_1)$ and $x_{t}(\xi_{2},u_{2},\Delta_2)$ of (\ref{TDS}) exist for any $t\geq 0$ and, furthermore, the following inequality holds:

\begin{eqnarray}
\label{ineqdeltaidss}
\begin{array}{rcl}
\left\Vert x_{t}(\xi_{1},u_{1},\Delta_1)-x_{t}(\xi_{2},u_{2},\Delta_2)\right\Vert_{\infty} & \leq &
\max\{\beta(\left\Vert \xi_{1}-\xi_{2}\right\Vert_{\infty}
,t),\\
& &
\gamma_{U}(\left\Vert (u_{1}-u_{2})|_{[-r,t-r[} \right\Vert
_{\infty}) 
\\
& &
+ \gamma_{D}(\Vert (\Delta_{1}-\Delta_{2})|_{[0,t[ } \Vert _{\infty} )\}.
\end{array}
\end{eqnarray}
\end{definition}

The above notion extends the one proposed in \cite{PolaSCL10} for time--delay systems with known and constant delays, to time--delay systems with unknown and time--varying delays. The main novelty relies on the last term in the right-hand side of (\ref{ineqdeltaidss}). The introduced stability notion regards solutions obtained by different time--delay signal realizations, and not only by different initial conditions and inputs. 
In general the inequality in (\ref{ineqdeltaidss}) is difficult to check
directly. We therefore provide hereafter a characterization of
$\delta$--IDSS, in terms of Lyapunov--Krasovskii functionals.

\begin{definition}\label{definitionLKfunctional}
Given a time--delay system $\Sigma=(X,\mathcal{X},\xi_{0},U,\mathcal{U},\mathcal{D},f)$, a locally
Lipschitz functional
$$
V:\mathbb{R}^{+}_{0}\times  C^{0}([-\Delta_{\max},0];\mathbb{R}^{n})\times
C^{0}([-\Delta_{\max},0];\mathbb{R}^{n})\to \mathbb{R}^{+}_{0}
$$
is said to be a $\delta$--IDSS Lyapunov--Krasovskii functional for
$\Sigma$ if there exist $\mathcal{ K}_{\infty}$ functions
$\alpha_1, \alpha_2$ and $\mathcal{K}$ functions $\alpha_3$,
$\rho$ such that:

\begin{itemize}

\item[(i)] for all $t\in \mathbb{R}_0^+$, $x_{1},x_{2}\in C^{0}([-\Delta_{\max},0];\mathbb{R}^{n})$
\[
\alpha_1(\Vert x_{1}(0)-x_{2}(0) \Vert)\le V(t,x_1,x_2)\le
\alpha_2(M_a(x_1-x_2)),
\]
where $M_{a}:C^{0}([-\Delta_{\max},0];\mathbb{R}^{n}) \to
\mathbb{R}_0^+$ is a continuous functional such that
\[
\underline\gamma_{a}(\Vert x(0)\Vert )\le M_{a}(x)\le \overline
\gamma_{a}(\Vert x\Vert_{\infty}),
\]
for all $x\in C^{0}([-\Delta_{\max},0];\mathbb{R}^{n})$, for some $\mathcal{K}_{\infty}$ functions $\underline \gamma_{a}$ and $\overline \gamma_{a}$;

\item[(ii)] by setting $L=\max\{M_1,\sup_{x,y\in X,u\in U} \Vert f(x,y,u)\Vert\}$, for all $x_{1},x_{2}\in \mathcal{X}$, $u_{1}$, $u_{2}\in U$, $d\in\mathcal{B}_{L(\Delta_{\max}-\Delta_{\min})}(0)\subset \mathbb{R}^{n}$
for which
\[
M_a(x_1-x_2)\geq \rho\left(\left \Vert (u_1-u_{2},d)\right \Vert\right),
\]
the following inequality
holds, almost everywhere in $t\in \mathbb{R}_0^+$:
\[
D^+V(t,x_1,x_2,u_1,u_2,d)\le -\alpha_3(M_a(x_1-x_2)),
\]
where $D^+V(t,x_1,x_2,u_1,u_2,d)$ is the derivative of the functional $V$ defined as
\begin{equation}
D^+V(t,x_1,x_{2},u_{1},u_{2},d)=
\limsup_{\theta\to 0^+} \frac {V\left (t+\theta,x^{\theta,\Delta(t)}_1,x^{\theta,\Delta(t)}_2\right )-V(t,x_1,x_2)}{\theta},\nonumber
\end{equation}
with

\begin{equation}
\begin{array}
{l}
x_1^{\theta,\Delta(t)}(s)=\left
\{
\begin{array}
{l}
x_1(s+\theta), \text{   if } s\in [-\Delta_{\max},-\theta[, \\
x_1(0)+(s+\theta)f(x_1(0),x_1(-\Delta(t)) +d,u_1), \text{   if } s\in [-\theta,0],
\end{array}
\right.\nonumber
\\
\\
x_2^{\theta,\Delta(t)}(s)=\left
\{
\begin{array}{l}
x_2(s+\theta), \text{   if } s\in [-\Delta_{\max},-\theta [, \\
x_2(0)+(s+\theta)f(x_2(0),x_2(-\Delta(t)),u_2), \text{   if } s\in [-\theta,0].
\end{array}
\right.\nonumber
\end{array}
\end{equation}

\end{itemize}
\label{KLF}
\end{definition}

Before providing a characterization of $\delta$--IDSS in terms of Lyapunov--Krasovskii functionals we need some preliminary technical results that we report hereafter.

\begin{lemma}\label{meanvaluetheorem}
Let $c\in\mathbb{R}$, $p\in\mathbb{N}$ and $g:[c,+\infty[\to
\mathbb{R}^p$ be a locally absolutely continuous function. Then,
for any $\xi_1, \xi_2 \in [c,+\infty[$ with $\xi_2\ge \xi_1$, the
following equality holds:
\begin{equation}
g(\xi_2)-g(\xi_1) =\left (  \int_0^1\left . \frac{dg(\xi)}{d\xi}\right |_{\xi= \xi_1+t (\xi_2-\xi_1)}dt\right ) (\xi_2-\xi_1).\nonumber
\end{equation}
\end{lemma}

\begin{lemma}\label{lemmad}
For any solution $x(t,\xi_0,u,\Delta)$ at time $t\geq
-\Delta_{\max}$ of system (\ref{TTDS}) and for any $a,b \in [-\Delta_{\max},+\infty[$ the following inequality holds:
\begin{equation}
\Vert x(b,\xi_0,u,\Delta)-x(a,\xi_0,u,\Delta)\Vert\le L \vert b-a\vert,
\end{equation}
where we recall $L=\max\{M_1,\sup_{x,y\in X,u\in U} \Vert f(x,y,u)\Vert\}$.
\end{lemma}
\begin{proof}
By Assumptions (A.0) and (A.1) and the definition of $X$, any solution lives in $X$ for any time  $t\geq -\Delta_{\max}$, and therefore the norm of the derivative of the solution in $[-\Delta_{\max},+\infty[$ is upper bounded by $\sup_{x,y\in X,u\in U} \Vert f(x,y,u) \Vert$ in $\mathbb{R}^+$ and by $M_1$ in $[-\Delta_{\max},0]$. Then, the result holds as a direct application of Lemma \ref{meanvaluetheorem}.
\end{proof}

\bigskip
We now have all the ingredients to present the following result.
\bigskip

\begin{theorem}
A time--delay system $\Sigma$ is $\delta$--IDSS if it admits a $\delta$--IDSS Lyapunov--Krasovskii functional.
\label{Th_KLF}
\end{theorem}
\begin{proof}
Let $\Delta_1$, $\Delta_2$ be two realizations of time--delays,
$\phi_1$, $\phi_2$ be two initial conditions, $u_1, u_2$ be two
inputs. Let $x_i(t)$, $i=1,2$, be the solutions of (\ref{TDS})
corresponding to $\phi_i$, $u_i$, $\Delta_i$, respectively. Consider:
\begin{eqnarray}&&
\dot x_1(t)=f(x_1(t),x_1(t-\Delta_2(t))+d(t),u_1(t-r)), \nonumber \\
&& d(t)=x_1(t-\Delta_1(t))-x_1(t-\Delta_2(t)),\nonumber \\ && \dot
x_2(t)=f(x_2(t),x_2(t-\Delta_2(t)),u_2(t-r)). \nonumber
\end{eqnarray}
By Lemma \ref{lemmad} it follows that
\begin{eqnarray}\label{upperboundd}
\Vert d(t)\Vert\le L\Vert \Delta_1(t)-\Delta_2(t) \Vert \le L(\Delta_{\max}-\Delta_{\min}).
\end{eqnarray}
By (\ref{upperboundd}), the proof is over if the existence of the
$\delta$--IDSS Lyapunov--Krasovskii functional guarantees that the system described
by
\begin{eqnarray}
\dot x(t)=f(x(t),x(t-\Delta(t))+d(t),u(t-r))\nonumber
\end{eqnarray}
is such that there exist a $\mathcal{KL}$ function $\bar \beta$,
${\mathcal K}$ functions $\bar \gamma_U$, $\bar \gamma_D$ such
that, for any realization of the time--delay $\Delta(t)$, initial
conditions $\phi_i$, inputs $u_i(t)$, $i=1,2$ and $d(t)$ (with
their upper bounds), the following inequality holds for the solutions,

\begin{equation}
\label{deltaISS}
\begin{array}
{rcl}
\Vert x_t(\phi_1,u_1,d)-x_t(\phi_2,u_2,0)\Vert & \le & 
\max \{\bar \beta (\Vert \phi_1-\phi_2\Vert_{\infty},t),\bar
\gamma_U\left (\left \Vert (u_1-u_2)_{[-r,t-r}\right
\Vert_{\infty}\right)
\\
& &
+\bar \gamma_D\left (\left \Vert d_{[0,t)}\right
\Vert_{\infty}\right)\}.
\end{array}
\end{equation}

We now show that the same lines of the proof used by Sontag for ISS (see, for instance, \cite{Sontag95}), used also for time--delay systems in \cite{PepeJiangSCL}, can be used here. Let $\phi_1,\phi_2$ be a pair of initial conditions. Let $u_1,u_2,d$ be such that $\left \Vert (u_1-u_2,d)
\right \Vert_{\infty} = v$. The set
\[
S=\{(t,\psi_1,\psi_2)\in \mathbb{R}^+\times
C^{0}([-\Delta_{\max},0];\mathbb{R}^{n})\times
C^{0}([-\Delta_{\max},0];\mathbb{R}^{n}): V(t,\psi_1,\psi_2)\le
\alpha_2\circ \rho(v) \}
\]
is such that: if $(t_0, x_{t_0}(\phi_1,u_1,d), x_{t_0}(\phi_2,u_2,0))\in S$ for some $t_0\in\mathbb{R}^{+}_{0}$, then
\[
(t,x_{t}(\phi_1,u_1,d),x_{t}(\phi_2,u_2,0))\in S,
\]
for all $t\ge t_0$. To prove this, by contradiction suppose that there exists a time $\bar t> t_0$ such that $V(\bar t,x_{\bar t}(\phi_1,u_1,d), x_{\bar t}(\phi_2,u_2,0))>\alpha_2\circ \rho (v)$. Then there exists a positive real $\epsilon$ and a time
$t_1>t_0$ such that
\[
V(t_{1},x_{t_1}(\phi_1,u_1,d), x_{t_1}(\phi_2,u_2,0))\ge \alpha_2\circ \rho (v)+\epsilon.
\]
Let $t_1$ be minimal for this given $\epsilon$. By continuity arguments, it follows that there exists $a\in\mathbb{R}^{+}$ such that for any $t\in ]t_1-a,t_1+a[$, $V(t,x_{t}(\phi_1,u_1,d), x_{ t}(\phi_2,u_2,0))>\alpha_2\circ \rho (v)$. By (ii), Theorem 5 in \cite{PepeTACABS} and Theorem 2 in \cite{PepeSAFEAUTOMATICA}, it follows that $V(t,x_{t}(\phi_1,u_1,d), x_{t}(\phi_2,u_2,0))$ is nonincreasing in $]t_1-a,t_1]$. Thus $t_1$ is not minimal and a contradiction holds. Now two cases can occur: 1) $(0,\phi_1,\phi_2)$ belongs to $S$; 2) $(0,\phi_1,\phi_2)$ does not belong to $S$. In case 1) by the forward invariant property of the set $S$, the following inequality holds, for any $t\geq 0$,
\[
\Vert x(t,\phi_1,u_1,d)-x(t,\phi_2,u_2,0)\Vert \le
\alpha_1^{-1}\circ \alpha_2\circ \rho(v).
\]
In case 2) let $t_{0}\in\mathbb{R}^{+}_0\cup \{+\infty\}$ be the maximal time (notice that it can be $+\infty$) such that, in the interval $[0,t_0[$, $(t, x_{t}(\phi_1,u_1,d ), x_{t}(\phi_2,u_2,0))$ does not belong to $S$. Then, in $[0,t_0[$, the inequality in (ii) holds, that, by Theorem 2 in \cite{PepeSAFEAUTOMATICA}, yields for $w(t)=V(t,x_t(\phi_1,u_1,d),$ $x_t(\phi_2,u_2,0))$, the inequality $D^+w(t)\le -\alpha_3\circ \alpha_2^{-1}(w(t))$ a.e. From this inequality, by (i), Theorem 5 in \cite{PepeTACABS}, the following inequality holds, for a suitable $\mathcal{KL}$ function $\tilde \beta$ (see Lemma 4.4 in \cite{SmoothConverse}), for $t\in [0,t_0[$,
\begin{eqnarray}
&& \Vert x(t,\phi_1,u_1,d)-x(t,\phi_2,u_2,0)\Vert\le 
 \alpha_1^{-1}\circ \tilde \beta (\alpha_2\circ \overline
\gamma_a(\Vert \phi_1-\phi_2 \Vert_{\infty}),t).
\end{eqnarray}
By the result concerning the forward invariant property of the set $S$, the following inequality holds, $\forall t\ge 0$:
\begin{equation}
\label{connormaeuclideasontag}
\begin{array}
{rcl}
\Vert
x(t,\phi_1,u_1,d)-x(t,\phi_2,u_2,0)\Vert & \le & 
\max \{\alpha_1^{-1}\circ \tilde \beta (\alpha_2\circ \overline
\gamma_a(\Vert \phi_1-\phi_2\Vert_{\infty}),t),\\
& & \alpha_1^{-1}\circ \alpha_2\circ \rho(v)\}.
\end{array}
\end{equation}

Therefore, in both cases 1) and 2), the inequality (\ref{connormaeuclideasontag}) holds for any $t\in\mathbb{R}^{+}$. From the inequality (\ref{connormaeuclideasontag}), one gets:
\begin{equation}
\begin{array}
{rcl}
\Vert x_t(\phi_1,u_1,d)-x_t(\phi_2,u_2,0)\Vert_{\infty} & \le & 
\max \{e^{-(t-\Delta_{\max})}\Vert \phi_1-\phi_2\Vert_{\infty}
\\
& & +  \alpha_1^{-1}\circ \tilde \beta (\alpha_2\circ \overline
\gamma_a(\Vert \phi_1-\phi_2\Vert_{\infty}),\\
& &
\max
\{0,t-\Delta_{\max}\}),
\alpha_1^{-1}\circ \alpha_2\circ \rho(v)\}
\end{array}
\end{equation}

By causality arguments (see \cite{Sontag95}), i.e. the
solutions $x_t(\phi_1,u_1,d)$ and $x_t(\phi_2,$ $u_2,0)$ are
independent of the inputs $u_1$, $u_2$ at time greater than $t-r$
(recall the input delay), and of the input $d$ at time greater
than $t$, it is possible to replace $v$ by $\Vert
((u_{1}-u_{2})|_{[-r,t-r[},d|_{[0,t[}) \Vert_{\infty}$. Thus, the
inequality in (\ref{deltaISS}) is proved, which concludes the
proof.
\end{proof}

\bigskip
\begin{remark}
Sufficient and necessary conditions for a time--delay system to be ISS, in terms of existence of ISS Lyapunov--Krasovskii functionals, can be found in \cite{PepeEJC08} and in \cite{cinesi}. For a time--delay system to be $\delta$--IDSS it is unclear at this stage whether the existence of a $\delta$--IDSS Lyapunov--Krasovskii functional is not only a sufficient condition but also a necessary one. 
When inputs are not involved (i.e. $u_1(t)=u_2(t)\equiv 0$, $\forall t\ge 0$), and the same realization of time--delay signal is considered (i.e. $\Delta_1(t)=\Delta_2(t)$, $\forall t\ge 0$), the notion of \mbox{$\delta$--IDSS} reduces to the well known \mbox{$\delta$--GAS} notion (see \cite{incrementalS}), extended to nonlinear time--delay systems. The \mbox{$\delta$--GAS} notion can be also defined as a GAS notion (see \cite{incrementalS}) with respect to a suitable diagonal set  by means of an extended system constructed as follows: the original system is replaced by two subsystems, each one described by the same dynamic functions of the original one. To this aim, the following result, that can be regarded as an extension of Lemma 2.3 in \cite{incrementalS}, is useful. Let the product space \mbox{$\mathbb{S}=C^0([-\Delta_{max},0]; \mathbb{R}^n)\times C^0([-\Delta_{max},0];\mathbb{R}^n)$} be endowed with the norm given, for $\psi=\left [\begin{array}{c} \psi_1\\ \psi_2\end{array}\right ]\in \mathbb{S}$, $\psi_i\in C^0([-\Delta_{max},0];\mathbb{R}^n)$, $i=1,2$, as $\Vert\psi\Vert_{\mathbb{S}}=\Vert \psi_1\Vert_{\infty}+\Vert \psi_2\Vert_{\infty}$. 
Let $\mathbb{H}$ be the subset (diagonal set) of $\mathbb{S}$ defined as $\left \{\left [\begin{array}{c} \psi\\ \psi\end{array}\right ], \psi \in C^0([-\Delta_{max},0]; \mathbb{R}^n)\right \}$. 
Then, for any $\psi_1, \psi_2 \in C^0([-\Delta_{max},0]; \mathbb{R}^n)$, the equality holds
\begin{equation}\label{angelilemma}
\Vert \psi_2-\psi_1\Vert_{\infty}=\inf_{\chi\in C^0([-\Delta_{max},0];\mathbb{R}^n)}\left \Vert\left [ \begin{array}{c}\psi_1-\chi \\ \psi_2-\chi\end{array} \right ]\right \Vert_{\mathbb{S}}.
\end{equation}
The proof of (\ref{angelilemma}) is as follows. 
Take any $\psi=\left [\begin{array}{c} \psi_1\\ \psi_2\end{array}\right ]\in \mathbb{S}$. The inequality holds
\begin{eqnarray} \label{primaineq}
\inf_{\chi\in C^0([-\Delta_{max},0];\mathbb{R}^n)}\left \Vert\psi-\left [ \begin{array}{c}\chi \\ \chi\end{array} \right ]\right \Vert_{\mathbb{S}}\le \left \Vert \psi-
\left [ \begin{array}{c}\psi_1 \\ \psi_1\end{array} \right ]
\right \Vert_{\mathbb{S}}=\Vert \psi_2-\psi_1\Vert_{\infty}.
\end{eqnarray}
On the other hand,
\begin {eqnarray}\label {secondaineq}&&
\Vert \psi_2-\psi_1\Vert_{\infty} \le \inf_{\chi\in  C^0([-\Delta_{max},0];\mathbb{R}^n)} (\Vert \psi_2-\chi\Vert_{\infty}+\Vert \psi_1-\chi \Vert_{\infty})=\\ \nonumber && \inf_{\chi\in C^0([-\Delta_{max},0];\mathbb{R}^n) } \left \Vert \psi-\left [\begin{array}{c} \chi\\ \chi\end{array}\right ]\right \Vert_{\mathbb{S}}.
\end{eqnarray}
From (\ref{primaineq}) and (\ref{secondaineq}) the equality (\ref{angelilemma}) follows. This result is in no contradiction with the result in 
Lemma 2.3 in \cite{incrementalS} (where a factor $\sqrt{2}/2$ is involved). Actually the same result shown in (\ref{angelilemma})
holds for the finite dimensional case, when the norm in $\mathbb{R}^{2n}=\mathbb{R}^n\times \mathbb{R}^n$ is chosen as the sum of the euclidean norms in $\mathbb{R}^n$.
Notice that the right-hand side of (\ref{angelilemma}) is a point distance of $\left [ \begin{array}{c}\psi_1\\ \psi_2\end{array} \right ]$ from the diagonal set $\mathbb{H}$. Therefore, the \mbox{$\delta$--GAS} of the original system is equivalent to \mbox{GAS} property of the above mentioned extended system, with respect to the diagonal set $\mathbb{H}$, also for time--delay systems.
By using this result, the converse Lyapunov--Krasovskii theorems available in the literature for the \mbox{GAS} property (see \cite{Karaf08a},\cite{KarJiaBook}), extended to the case of stability with respect to the diagonal set, can be used directly to provide converse Lyapunov--Krasovskii theorems for the \mbox{$\delta$--GAS} notion of time--delay systems. Whether the results in \cite{PepeEJC08} could be used also for the \mbox{$\delta$--IDSS}, in order to provide converse Lyapunov--Krasovskii theorems, is an interesting and not easy topic of future investigations. The main issue concerns dealing with different time--delay signal realizations, which to our knowledge has never been considered in the literature so far. 
\end{remark}

Quadratic Lyapunov--Krasovskii functionals \cite{TDSbook} are popularly used for the stability check of time--delay systems. They are an important set of functionals that can be useful also for checking the $\delta$--IDSS property, as it will be shown in Section 7 through an illustrative example.

\section{Symbolic Models and Approximate Equivalence Notions}

We use the class of alternating transition systems \cite{Alternating} as abstract mathematical models of time--delay systems.

\begin{definition}
\label{Def_TS}
An (alternating) transition system is a tuple:
\[
T=(Q,q_{0},L,\rTo,O,H),
\]
consisting of:

\begin{itemize}
\item A set of states $Q$;

\item An initial state $q_{0}\in Q$;

\item A set of labels $L=A\times B$, where:

\begin{itemize}
\item $A$ is the set of control labels;

\item $B$ is the set of disturbance labels;

\end{itemize}

\item A transition relation $\rTo\subseteq Q\times L\times Q$;

\item An output set $O$;

\item An output function $H:Q\rightarrow O$.
\end{itemize}

A transition system $T$ is said to be:

\begin{itemize}
\item \textit{metric}, if the output set $O$ is equipped with a metric $\mathbf{d}:O\times O\rightarrow\mathbb{R}_{0}^{+}$;

\item \textit{countable}, if $Q$ and $L$ are countable sets;

\item \textit{finite}, if $Q$ and $L$ are finite sets.
\end{itemize}
\end{definition}

We follow standard practice and denote by $q \rTo^{a,b} p$, a transition from $q$ to $p$ labeled by $a$ and $b$. Transition systems capture dynamics through the transition relation. For any states $q,p\in Q$, the symbol $q \rTo^{a,b} p$ simply means that it is possible to evolve or jump from state $q$ to state $p$ under the action labeled by $a$ and $b$. \\
We consider simulation and bisimulation relations \cite{Milner,Park} as formal equivalence notions to relate time--delay systems and their symbolic models. Intuitively, a bisimulation relation between a pair of transition systems $T_{1}$ and $T_{2}$ is a relation between the corresponding state sets explaining how a state trajectory $r_{1}$ of $T_{1}$ can be transformed into a state trajectory $r_{2}$ of $T_{2}$, and vice versa. While typical bisimulation relations require that $r_{1}$ and $r_{2}$ have the same output run, i.e. $H_{1}(r_{1}) = H_{2}(r_{2})$, the notions of approximate bisimulation relations, that we now introduce, relax this condition and require that $H_{1}(r_{1})$ is simply close to $H_{2}(r_{2})$ where closeness is measured with respect to a metric on the output set. A simulation relation is a one--sided version of a bisimulation relation. 

\begin{definition}
\label{ASR1}
\cite{AB-TAC07} 
Given two metric transition systems $T_{1}=(Q_{1},q_{1}^{0},A_{1}\times B_{1},\rTo_{1},$ $O,H_{1})$ and $T_{2}=(Q_{2},q_{2}^{0},A_{2}\times B_{2},\rTo_{2},O,H_{2})$ with the same output set $O$ and the same metric $\mathbf{d}$
and given a precision $\varepsilon\in \mathbb{R}^{+}_{0}$, a relation $R\subseteq Q_{1}\times Q_{2}$ is said to be an $\varepsilon$--approximate simulation relation from $T_{1}$ to $T_{2}$ if for any $(q_{1},q_{2})\in R$:
\begin{itemize}
\item[(i)]$\mathbf{d}(H_{1}(q_{1}),H_{2}(q_{2}))\leq\varepsilon$;
\item[(ii)] $\forall (a_{1},b_{1})\in A_{1}\times B_{1}$, $\exists (a_{2},b_{2})\in A_{2}\times B_{2}$ such that $q_{1} \rTo^{a_{1},b_{1}}_{1} p_{1}$ and $q_{2} \rTo^{a_{2},b_{2}}_{2} p_{2}$ with $(p_{1},p_{2})\in R$.
\end{itemize}

Moreover, $T_{1}$ is said to be $\varepsilon$--approximately simulated from $T_{2}$, denoted $T_{1}\preceq_{\varepsilon} T_{2}$, if there exists an $\varepsilon$--approximate simulation relation $R$ from $T_{1}$ to $T_{2}$ such that $(q_{1}^{0},q_{2}^{0})\in R$. 
%
%
Given a precision $\varepsilon\in \mathbb{R}^{+}_{0}$, a relation $R\subseteq Q_{1}\times Q_{2}$ is said to be an $\varepsilon$--approximate bisimulation relation between $T_{1}$ and $T_{2}$ if:
\begin{itemize}
\item[(i)] $R$ is an $\varepsilon$--approximate simulation relation from $T_{1}$ to $T_{2}$;
\item[(ii)] $R^{-1}$ is an $\varepsilon$--approximate simulation relation from $T_{2}$ to $T_{1}$.
\end{itemize}

Moreover, $T_{1}$ is said to be $\varepsilon$--bisimilar to $T_{2}$, denoted $T_{1}\cong_{\varepsilon} T_{2}$, if there exists an $\varepsilon$--approximate bisimulation relation $R$ between $T_{1}$
and $T_{2}$ such that $(q_{1}^{0},q_{2}^{0})\in R$.
\end{definition}

\bigskip
The above notions have been further extended in \cite{PolaSIAM2009} to the ones of alternating approximate simulation and bisimulation, that appropriately capture the different role played by control labels and disturbance labels.

\begin{definition}
\label{Alt_ASR1}
Given two metric transition systems $T_{1}=(Q_{1},q_{1}^{0},A_{1}\times B_{1},\rTo_{1},$ $O,H_{1})$ and $T_{2}=(Q_{2},q_{2}^{0},A_{2}\times B_{2},\rTo_{2},O,H_{2})$ with the same output set $O$ and the same metric $\mathbf{d}$
and given a precision $\varepsilon\in \mathbb{R}^{+}_{0}$, a relation $R\subseteq Q_{1}\times Q_{2}$ is said to be an alternating $\varepsilon$--approximate ($A\varepsilon A$)
simulation relation from $T_{1}$ to $T_{2}$ if for any $(q_{1},q_{2})\in R$:
\begin{itemize}
\item[(i)]$\mathbf{d}(H_{1}(q_{1}),H_{2}(q_{2}))\leq\varepsilon$;
\item[(ii)] $\forall a_{1}$ $\in A_{1}$ $\exists a_{2}\in A_{2}$ $\forall b_{2}\in B_{2}$ $\exists b_{1}\in B_{1}$ such that $q_{1} \rTo^{a_{1},b_{1}}_{1} p_{1}$ and $q_{2} \rTo^{a_{2},b_{2}}_{2} p_{2}$ with $(p_{1},p_{2})\in R$.
\end{itemize}

Moreover, $T_{1}$ is said to be $A\varepsilon A$ simulated from $T_{2}$, denoted $T_{1}\preceq_{\varepsilon}^{\alt} T_{2}$, if there exists an $A\varepsilon A$ simulation relation $R$ from $T_{1}$ to $T_{2}$ such that $(q_{1}^{0},q_{2}^{0})\in R$. 
%
Given a precision $\varepsilon\in \mathbb{R}^{+}_{0}$, a relation $R\subseteq Q_{1}\times Q_{2}$ is said to be an alternating $\varepsilon$--approximate ($A\varepsilon A$)
bisimulation relation between $T_{1}$ and $T_{2}$ if:
\begin{itemize}
\item[(i)] $R$ is an alternating $\varepsilon$--approximate simulation relation from $T_{1}$ to $T_{2}$;
\item[(ii)] $R^{-1}$ is an alternating $\varepsilon$--approximate simulation relation from $T_{2}$ to $T_{1}$.
\end{itemize}

Moreover, $T_{1}$ is said to be $A\varepsilon A$ bisimilar to $T_{2}$, denoted $T_{1}\cong_{\varepsilon}^{\alt} T_{2}$, if there exists an $A\varepsilon A$ bisimulation relation $R$ between $T_{1}$
and $T_{2}$ such that $(q_{1}^{0},q_{2}^{0})\in R$.
\end{definition}

\bigskip
When $\varepsilon=0$, the notion of bisimulation in the above definition coincides with the $2$--players version of the definition proposed in \cite{Alternating}. 

\section{Existence of Symbolic Models}
\label{sec6}
\subsection{Alternating Approximately Bisimilar Symbolic Models}\label{sec5.2}
In this section we propose symbolic models that approximate $\delta$--IDSS time--delay systems in the sense of A$\varepsilon$A bisimulation.
These symbolic models are obtained by a first time discretization of the time--delay system and a subsequent quantization of the state and input variables, and of the delay signals. \\
In many concrete applications controllers are implemented through digital devices and this motivates our interest for
control inputs that are piecewise constant. In the following we refer to time--delay systems with digital controllers as \textit{digital time--delay systems}. We suppose that the set $U$ of input values contains the origin and that the class of control inputs is
\begin{equation}
\mathcal{U}_{\tau}:=
\{
u:[-r,-r+\tau ] \rightarrow U \,|\, \exists v\in\mathcal{U} \text{ s.t. } u=v|_{[-r,-r+\tau]}
\},
\label{Utau}
\end{equation}

where \mbox{$\tau\in\mathbb{R}^{+}$} can be thought of as a sampling time parameter in the digital controller. For further purposes, define:
\[
\mathcal{D}_{\tau}:=
\{
\Delta:[0,\tau ] \rightarrow [\Delta_{\min},\Delta_{\max}] \,|\, \exists d\in\mathcal{D} \text{ s.t. } \Delta=d|_{[0,\tau]}
\}.
\]
Given $k\in\mathbb{N}$ we denote by $\mathcal{U}_{k,\tau}$ and by $\mathcal{D}_{k,\tau}$ the class of control inputs and time--delay signals obtained by the concatenation of $k$ control inputs in $\mathcal{U}_{\tau}$ and respectively, of $k$ time--delay signals in $\mathcal{D}_{\tau}$. Given a digital time--delay system $\Sigma$ define the transition system
\[
T_{\tau}(\Sigma):=(Q_{1},q_{1}^{0},L_{1},\rTo_{1},O_{1},H_{1}),
\]
where:
\begin{itemize}
\item $Q_{1}=\mathcal{X}$;
\item $q_{1}^{0}=\xi_{0}$;
\item $L_{1}=A_{1}\times B_{1}$, where:
\begin{itemize}
\item $A_{1}=\{a_{1}\in \mathcal{U}_{\tau}\,\,\vert\,\,x_\tau (x,a_{1},b_{1})$ is defined for any $x\in\mathcal{X} \text{ and } b_{1}\in B_{1}\}$;
\item $B_{1}=\{b_{1}\in \mathcal{D}_{\tau}\,\,\vert\,\,x_\tau (x,a_{1},b_{1})$ is defined for any $x\in\mathcal{X} \text{ and } a_{1}\in A_{1}\}$;
\end{itemize}
\item $q \rTo_{1}^{a_{1},b_{1}} p$, if
$x_{\tau}(q,a_{1},b_{1})=p$;
\item $O_{1}=\mathcal{X}$;
\item $H_{1}=1_{X}$.
\end{itemize}
Transition system $T_{\tau}(\Sigma)$ can be thought of as a time discretization of $\Sigma$ and it is metric when we regard $O_{1}=\mathcal{X}$ as being equipped with the metric \mbox{$\mathbf{d}(p,q)=\Vert p-q\Vert_{\infty}$}. Note that $T_{\tau}(\Sigma)$ is not symbolic, since the set of states $Q_{1}$ is a functional space. The construction of symbolic models for digital time--delay systems relies upon approximations of the set of reachable states, of
the space of time--delay signals and of the space of input signals. Let $R_{\tau}(\Sigma)\subseteq \mathcal{X}$ be the set of reachable states of $\Sigma$ at times $t=0,\tau,...,k\tau,...$, i.e. the collection of all states $x\in\mathcal{X}$ for which there exist $k\in \mathbb{N}$, a control input $u\in\mathcal{U}_{k,\tau}$ and a time--delay signal $\Delta\in \mathcal{D}_{k,\tau}$, so that $x=x_{k\tau}(\xi_{0},u,\Delta)$. The sets $R_{\tau}(\Sigma)$, $A_{1}$ and $B_{1}$ are functional spaces and are therefore needed to be approximated, in the sense of the following definition.

\begin{definition} \cite{PolaSCL10}
\label{CountApprox}
Consider a functional space $\mathcal{Y}\subseteq C^{0}(I,Y)$ with $Y\subseteq\mathbb{R}^{n}$, \mbox{$I=[i_{1},i_{2}]$}, $i_{1},i_{2}\in \mathbb{R}$, $i_{1}<i_{2}$.
A map \mbox{$\mathcal{A}:\mathbb{R}^{+}\rightarrow 2^{C^{0}(I,Y)}$} is a \textit{countable approximation} of $\mathcal{Y}$ if for any desired precision $\lambda\in\mathbb{R}^{+}$:
\begin{itemize}
\item [(i)] $\mathcal{A}(\lambda)$ is a countable set;
\item [(ii)] $\forall y\in \mathcal{Y}$, $\exists z\in \mathcal{A}(\lambda)$ so that $\Vert y - z \Vert_{\infty} \leq \lambda$;
\item [(iii)] $\forall z\in \mathcal{A}(\lambda)$,  $\exists y\in \mathcal{Y}$ so that $\Vert y - z \Vert_{\infty} \leq \lambda$.
\end{itemize}
\end{definition}
A countable approximation $\mathcal{A}_{U}$ of $\mathcal{U}_{\tau}$ can be easily obtained by defining for any $\lambda_{U}\in\mathbb{R}^{+}$,
\begin{eqnarray}&&
\mathcal{A}_{U}(\lambda_{U})=\{u\in\mathcal{U}_{\tau}:
u(t)=u(-r)\in [U]_{2\lambda_{U}}, 
t\in\lbrack -r,-r+\tau]\}, \label{UtauQ}
\end{eqnarray}
where $[U]_{2\lambda_{U}}$ is defined as in (\ref{grid}). 
By comparing $\mathcal{U}_{\tau}$ in (\ref{Utau}) and $\mathcal{A}_{U}(\lambda_{U})$ in (\ref{UtauQ}) it is readily seen that $\mathcal{A}_{U}(\lambda_{U})\subset \mathcal{U}_{\tau}$ for any $\lambda_{U}\in\mathbb{R}^{+}$.
Since the set $U$ contains the origin, the set $\mathcal{A}_{U}(\lambda_{U})$ is nonempty.  
Suppose that the time--delay system $\Sigma$ is $\delta$--IDSS and let $\beta$, $\gamma_{U}$ and $\gamma_{D}$ be a $\mathcal{KL}$ and $\mathcal{K}$ functions satisfying the inequality in (\ref{ineqdeltaidss}).
We can now define a countable transition system that approximates $T_{\tau}(\Sigma)$.
Given a time quantization $\tau\in\mathbb{R}^{+}$, a state space quantization $\lambda_{X}\in\mathbb{R}^{+}$, an input space quantization $\lambda_{U}\in\mathbb{R}^{+}$ and  a delay quantization $\lambda_{D}\in\mathbb{R}^{+}$, consider the transition system:
\begin{equation}
T_{\tau,\lambda_{X},\lambda_{U},\lambda_{D}}(\Sigma):=(Q_{2},q_{2}^{0},L_{2},\rTo_{2},O_{2},H_{2}),
\label{sm}
\end{equation}
where:
\begin{itemize}
\item $Q_{2}=\mathcal{A}_{X}(\lambda_{X})$;
\item $q_{2}^{0}\in Q_{2}$ such that $\Vert \xi_{0} - q_{2}^{0} \Vert_\infty \leq \lambda_{X}$;
\item $L_{2}=A_{2}\times B_{2}$, where $A_{2}=\mathcal{A}_\mathcal{U}(\lambda_{U})$ and $B_{2}=\mathcal{A}_{D}(\lambda_{D})$;
\item $q \rTo^{a_{2},b_{2}}_{2} p$, if
\[
\left\Vert p-x_{\tau}(q,a_{2},b_{2})\right\Vert_{\infty}\leq \max\{\beta(\lambda_{X},\tau),\gamma_{D}(\lambda_{D})\}+\lambda_{X};
\]
\item $O_{2}=\mathcal{X}$;
\item $H_{2}=\imath : Q_{2} \hookrightarrow O_{2}$.
\end{itemize}
We can now give the following result:

\begin{theorem}
Consider a digital time--delay system $\Sigma
$ and any desired precision $\varepsilon\in\mathbb{R}^{+}$. Suppose that $\Sigma$ is $\delta$--IDSS and choose $\tau\in\mathbb{R}^{+}$ so that $\beta(\varepsilon,\tau)<\varepsilon$. Moreover suppose that there exist countable approximations $\mathcal{A}_{X}$ and $\mathcal{A}_{D}$ of $R_{\tau}(\Sigma)$ and $\mathcal{D}_{\tau}$, respectively. Then, for any $\lambda_{X}\in\mathbb{R}^{+}$, $\lambda_{U}\in\mathbb{R}^{+}$ and $\lambda_{D}\in\mathbb{R}^{+}$ satisfying the following inequality:
\begin{eqnarray}
\max\{\beta(\varepsilon,\tau),\gamma_{U}(\lambda_{U})+ \gamma_{D}(\lambda_{D})\}+
\max\{\beta(\lambda_{X},\tau),\gamma_{D}(\lambda_{D})\}+\lambda_{X}
\leq \varepsilon,
\label{cond2}
\end{eqnarray}
transition systems $T_{\tau,\lambda_{X},\lambda_{U},\lambda_{D}}(\Sigma)$ and $T_{\tau}(\Sigma)$ are A$\varepsilon$A bisimilar, i.e.
$T_{\tau,\lambda_{X},\lambda_{U},\lambda_{D}}(\Sigma)$ $\cong_{\varepsilon}^{\alt} T_{\tau}(\Sigma)$.
\label{ThMain}
\end{theorem}

Before giving the proof of the above result we point out that:

\begin{lemma}
For any given precision $\varepsilon\in\mathbb{R}^{+}$, there exists a choice of quantization parameters $\tau,\lambda_{X},\lambda_{U},\lambda_{D}\in\mathbb{R}^{+}$ so that the inequality in (\ref{cond2}) holds.
\end{lemma}

\begin{proof}
Pick $\lambda_{X}\leq \varepsilon/3$. Since $\gamma_{U}$ and $\gamma_{D}$ are $\mathcal{K}$ functions, there exists a choice of $\lambda_{U}$ and $\lambda_{D}$ so that $\gamma_{U}(\lambda_{U})\leq \varepsilon/6$ and $\gamma_{D}(\lambda_{D})\leq \varepsilon/6$. Since $\beta$ is a $\mathcal{KL}$ function there exists $\tau$ so that $\beta(\varepsilon,\tau)\leq \varepsilon/3$. By construction, since $\lambda_{X}< \varepsilon$, then $\max\{\beta(\lambda_{X},\tau),\gamma_{D}(\lambda_{D})\}< \max\{\beta(\varepsilon,\tau),\gamma_{U}(\lambda_{U})+ \gamma_{D}(\lambda_{D})\}$, from which:
\begin{eqnarray}
&& \max\{\beta(\varepsilon,\tau),\gamma_{U}(\lambda_{U})+ \gamma_{D}(\lambda_{D})\}+
\max\{\beta(\lambda_{X},\tau),\gamma_{D}(\lambda_{D})\}+\lambda_{X}
\leq \nonumber \\
&& 2\max\{\beta(\varepsilon,\tau),\gamma_{U}(\lambda_{U})+ \gamma_{D}(\lambda_{D})\}+\lambda_{X}
\leq 
2\max\{\varepsilon/3,\varepsilon/6+\varepsilon/6\}+\varepsilon/3=\varepsilon, \nonumber
\end{eqnarray}
which concludes the proof.
\end{proof}

We can now give the proof of Theorem \ref{ThMain}. 

\begin{proof}
Consider the relation \mbox{$R\subseteq Q_{1}\times Q_{2}$} defined by $(x,q)\in R$ if and only if \mbox{$\Vert H_{1}(x)-H_{2}(q)\Vert_{\infty}\leq \varepsilon$}. We now show that $R$ is a A$\varepsilon$A simulation relation from $T_{\tau}(\Sigma)$ to  $T_{\tau,\lambda_{X},\lambda_{U},\lambda_{D}}(\Sigma)$. Consider any $(x,q)\in R$. Condition (i) in Definition \ref{Alt_ASR1} is satisfied by the definition of $R$. Let us now show that condition (ii) in Definition \ref{Alt_ASR1} holds.
Consider any $a_{1}\in A_{1}$. By definition of $A_{2}$ there exists $a_{2}\in A_{2}$ so that:
\begin{equation}
\left\Vert a_{1} - a_{2}\right\Vert_{\infty} \leq \lambda_{U}.
\label{no3}
\end{equation}
Consider any $b_{2}\in B_{2}$. By construction of $T_{\tau,\lambda_{X},\lambda_{U},\lambda_{D}}(\Sigma)$ there exists $b_{1}\in B_{1}$ so that:
\begin{equation}
\Vert b_{1}-b_{2} \Vert_{\infty}\leq \lambda_{D}.
\label{condelay}
\end{equation}

By definition of $Q_{2}$ there exists $\bar{x}\in R_{\tau}(\Sigma)$ so that:
\begin{equation}
\Vert \bar{x}-q \Vert_{\infty}\leq \lambda_{X}.
\label{cond111}
\end{equation}
Set $\bar{y}=x_{\tau}(\bar{x},a_{2},b_{1})\in R_{\tau}(\Sigma)$. By definition of $Q_{2}$ there exists $p\in Q_{2}$ so that:
\begin{equation}
\Vert \bar{y} - p \Vert_{\infty}\leq \lambda_{X}.
\label{cond222}
\end{equation}
Set $z=x_{\tau}(q,a_{2},b_{2})$. Note that since $a_{2}\in A_{2}\subseteq \mathcal{U}_{\tau}$, function $z$ is well defined. By the $\delta$--IDSS assumption and the inequalities in (\ref{cond111}) and (\ref{cond222}), the following chain of inequalities holds:
\begin{equation}
\begin{array}{rcl}
\label{no777}
\Vert z-p\Vert_{\infty} & = & \Vert z-\bar{y}+\bar{y}-p\Vert_{\infty}\\
& \leq & \Vert z-\bar{y}\Vert_{\infty}+\Vert \bar{y}-p \Vert_{\infty}\\
&\leq &  \max\{\beta(\Vert q-\bar{x}\Vert_{\infty},\tau),\gamma_{U}(\Vert a_{2}-a_{2}\Vert_{\infty})
+
\gamma_{D}(\Vert b_{1}-b_{2}\Vert_{\infty})\}+
\lambda_{X} \\
& \leq & \max\{\beta(\lambda_{X},\tau),\gamma_{U}
(0)+
\gamma_{D}(\lambda_{D})\}+\lambda_{X}.
\end{array}
\end{equation}
By the above inequality it is clear that $q \rTo^{a_{2},b_{2}}_{2} p$ in $T_{\tau,\lambda_{X},\lambda_{U},\lambda_{D}}(\Sigma)$.
Consider $x\rTo_{1}^{a_{1},b_{1}} y$. Since $\Sigma$ is $\delta$--IDSS
and by (\ref{cond2}), (\ref{no3}), (\ref{condelay}) and (\ref{no777}), the following chain of
inequalities holds:
\begin{equation}
\begin{array}{rcl}
\Vert y-p\Vert_{\infty} & \leq &
 \Vert y-z\Vert_{\infty}+\Vert z-p\Vert_{\infty}\\
 & \leq &
 \max\{\beta(\Vert x-q \Vert_{\infty},\tau),\gamma_{U}(\Vert a_{1}-a_{2} \Vert)
+
\gamma_{D}(\Vert b_{1}-b_{2} \Vert_{\infty})\} \nonumber\\
& &
+
\max\{
\beta(\lambda_{X},\tau),
\gamma_{D}(\lambda_{D})\}+\lambda_{X}\\
 & \leq &
 \max\{\beta(\varepsilon,\tau),\gamma_{U}(\lambda_{U})+\gamma_{D}(\lambda_{D})\} +
 \max\{\beta(\lambda_{X},\tau),
\gamma_{D}(\lambda_{D})\}+ \lambda_{X}\leq \varepsilon.
\label{nono7}
\end{array}
\end{equation}

Hence $(y,p)\in R$ and condition (ii) in Definition \ref{Alt_ASR1} holds. By the inequality in (\ref{cond2}) and the definition of $q_{2}^{0}$, \mbox{$\Vert \xi_{0}- q_{2}^{0} \Vert \leq \lambda_{X} \leq \varepsilon$} and hence, transition system
$T_{\tau}(\Sigma)$ is A$\varepsilon$A simulated by $T_{\tau,\lambda_{X},\lambda_{U},\lambda_{D}}(\Sigma)$. By using similar arguments it is possible to show that $R^{-1}$ is an A$\varepsilon$A simulation relation from $T_{\tau,\lambda_{X},\lambda_{U},\lambda_{D}}(\Sigma)$ to $T_{\tau}(\Sigma)$. Hence, the result follows.
\end{proof}

\subsection{Spline-based countable approximation of functional spaces}

The result presented in the previous section assumes existence of countable approximations of functional spaces of time--delay systems. By following \cite{PolaSCL10} in this section we present an approach to approximate these functional spaces which is based on spline analysis \cite{SplineBook}. Spline based approximation schemes have been extensively used in the literature of time--delay systems (see e.g. \cite{GermaniSIAM00} and the references therein).
Let us consider the space $\mathcal{Y}\subseteq C^{0}(I,Y)$ with $Y\subseteq\mathbb{R}^{n}$, $I=[i_{1},i_{2}]$, $i_{1},i_{2}\in\mathbb{R}$ and $i_{1}<i_{2}$.
Given $N\in \mathbb{N}$ consider the following functions (see \cite{SplineBook}):
\begin{equation}
\begin{array}{llll}
s_{0}(t)=
\left\{
\begin{array}
[c]{lll}
1-(t-i_{1})/h, & t\in [i_{1},i_{1}+h],&\\
0, & \textit{otherwise,} &
\end{array}
\right.
& & &  \\ \\
s_{i}(t)=
\left\{
\begin{array}
[c]{lll}
1-i+(t-i_{1})/h, & t\in [i_{1}+(i-1)h,i_{1}+ih], & \\
1+i-(t-i_{1})/h, & t\in [i_{1}+ih,i_{1}+(i+1)h], & \\
0, & \textit{otherwise,} 
\end{array}
\right.\\
\hspace{6mm} i=1,2,...,N; \\
\\
s_{N+1}(t)=
\left\{
\begin{array}
[c]{lll}
1+(t-i_{2})/h, &  t\in [i_{2}-h,i_{2}], & \\
0, & \textit{otherwise,} &
\end{array}
\right.
 & &  \\
\end{array}
\label{spline}
\end{equation}

where $h=(i_{2}-i_{1})/(N+1)$. Functions $s_{i}$ called \textit{splines}, are used to approximate $\mathcal{Y}$. 
Given any $N\in\mathbb{N}$, $\theta,M\in \mathbb{R}^{+}$ let be\footnote{The real $M$ is a parameter associated with $\mathcal{Y}$ and its role will become clear in the subsequent developments.}:
\begin{equation}
\Lambda(N,\theta,M):=h^{2} M/8+(N+2)\theta,
\label{lambda}
\end{equation}
with $h=(i_{2}-i_{1})/(N+1)$. Function $\Lambda$ will be shown to be an upper bound to the error associated with the approximation scheme that we propose. It is readily seen that for any $\lambda\in\mathbb{R}^{+}$ and any $M\in\mathbb{R}^{+}$ there always exist $N\in\mathbb{N}$ and $\theta\in\mathbb{R}^{+}$ so that $\Lambda(N,\theta,M)\leq \lambda$. Let $N_{\lambda,M}$ and $\theta_{\lambda,M}$ be such that $\Lambda(N_{\lambda,M},\theta_{\lambda,M},M)\leq \lambda$. For any $\lambda\in\mathbb{R}^{+}$ and $M\in\mathbb{R}^{+}$, define the operator $\psi_{\lambda,M}: \mathcal{Y} \rightarrow C^{0}(I;Y)$ that associates to any function $y\in\mathcal{Y}$ the function:
\begin{equation}
\psi_{\lambda,M}(y)(t):=\sum_{i=0}^{N_{\lambda,M}+1} \tilde{y}_{i}s_{i}(t), \hspace{5mm} t\in I,
\label{refpsi}
\end{equation}
where $\tilde{y}_{i}\in [Y]_{2\theta_{\lambda,M}}$ and $\Vert \tilde{y}_{i} - y(i_{1}+ih) \Vert \leq \theta_{\lambda,M}$, for any \mbox{$i=0,1,...,N_{\lambda,M}+1$}. Note that the operator $\psi_{\lambda,M}$ is not uniquely defined. For any given $M\in\mathbb{R}^{+}$ and any given precision $\lambda\in\mathbb{R}^{+}$ define:
\begin{equation}
\mathcal{A}_{\mathcal{Y},M}(\lambda):=\psi_{\lambda,M}(\mathcal{Y}).
\label{Ay}
\end{equation}

\begin{lemma}
\cite{PolaSCL10}
Suppose that $\mathcal{Y}\subseteq PC^{2}(I;Y)$ and there exists $M\in\mathbb{R}^{+}$ so that $\Vert D^{2}\,y\Vert_{\infty}\leq M$ for any $y\in \mathcal{Y}$. Then $\mathcal{A}_{\mathcal{Y},M}$ as defined in (\ref{Ay}), is a countable approximation of $\mathcal{Y}$.
\label{prop1}
\end{lemma}

The above result is useful for approximating the functional space $\mathcal{D}_{\tau}$ of time--delay signals, as shown hereafter:

\begin{proposition}
Consider the functional space $\mathcal{D}_{\tau}$ and $M_{D}\in\mathbb{R}^{+}$. Suppose that:
\begin{itemize}
\item[(A.2)] $\mathcal{D}_{\tau}=PC^{2}([0,\tau];[\Delta_{\min},\Delta_{\max}])$;
\item[(A.3)] for any $\Delta\in \mathcal{D}_{\tau}$, $\Vert D^{2} \Delta \Vert_{\infty}\leq M_{D}$.
\end{itemize}

Then the set $\mathcal{A}_{D}$ defined for any $\lambda_{D}\in\mathbb{R}^{+}$ by:
\begin{equation}
\mathcal{A}_{D}(\lambda_{D})=\psi_{\lambda_{D},M_{D}}(\mathcal{D}_{\tau}),
\label{ApproxDelay}
\end{equation}
with $\psi_{\lambda_{D},M_{D}}$ as in (\ref{refpsi}), is a countable approximation of $\mathcal{D}_{\tau}$.
\label{DelayApprox}
\end{proposition}

We now proceed with a further step towards the construction of countable approximations of $R_{\tau}(\Sigma)$.
Consider a digital time--delay system $\Sigma=(X,\mathcal{X},\xi_{0},U,\mathcal{U}_{\tau},\mathcal{D}_{\tau},f)$ and suppose that:
\begin{itemize}
\item[(A.4)] $\Sigma$ is $\delta$--IDSS;
\item[(A.5)] Function $f$ is Fr\'echet differentiable in ${\mathbb R}^n\times {\mathbb R}^n\times {\mathbb R}^m$;
\item[(A.6)] The Fr\'echet differential $J(x,y,u)$ of $f$ is continuous and bounded on bounded subsets of ${\mathbb R}^n\times {\mathbb R}^n\times {\mathbb R}^m$.
\end{itemize}
Under the above assumptions the following bounds are well defined:
\begin{equation}
\begin{array}
{cc}
B_{J}=\sup_{(x,y,u)\in X\times X\times U }\Vert J(x,y,u)\Vert,&
M_{X}=(2B_{X}+B_{U})(1+d_{\min})\kappa B_{J},
\label{cost1}
\end{array}
\end{equation}
where $\kappa$ is the Lipschitz constant of function $f$ in the
bounded set $X\times X\times U$ and $\Vert J(x,y,u)\Vert$ denotes
the norm of the operator $J(x,y,u): \mathbb{R}^{n}\times
\mathbb{R}^{n}\times \mathbb{R}^{m}\rightarrow \mathbb{R}^{n}$. We
can now give the following technical lemma that is instrumental
to prove the main result of this section.

\begin{lemma}
Consider a digital time--delay system $\Sigma$ satisfying assumptions (A.0), (A.1), (A.4--6), and
\begin{itemize}
\item[(A.7)] the following conditions:
\begin{itemize}
\item[(A.7.1)] $\xi_0\in C^0([-\Delta_{\max},0];X)\cap PC^2([-\Delta_{\max},0];\mathbb{R}^{n})$;
\item[(A.7.2)] $\left \Vert D^2\xi_0\right \Vert_{\infty}\le M_{X}$;
\item[(A.7.3)] $\beta_{ISS}(B_{X}^{0},\tau)+\gamma_{ISS}(B_{U})\leq B_{X}^{0}$;
\item[(A.7.4)] $\tau>2\Delta_{\max}$.
\end{itemize}

\end{itemize}

Then, for any $x_{k\tau}\in R_{\tau}(\Sigma)$ with $k\in\mathbb{N}$ the following hold:

\begin{itemize}
\item $x_{k\tau}\in C^0([-\Delta_{\max},0];X)\cap PC^2([-\Delta_{\max},0];\mathbb{R}^n)$;
\item $\Vert x_{k\tau} \Vert_{\infty}\leq B_{X}^{0}$;
\item $\left \Vert D^{2} x_{k\tau}\right \Vert_{\infty} \le M_{X}$.
\end{itemize}
\label{prop22}
\end{lemma}
\begin{proof} It is sufficient to show that $x_{\tau}$ satisfies the same properties of $\xi_0$,
i.e. conditions (A.7.1), (A.7.2) hold with $\xi_0$ replaced by
$x_{\tau}$ and $\Vert x_{\tau}\Vert_{\infty}\le B_X^0$. First note
that the function $t\to \dot{x}(t)$, $t\in [0,\tau]$, is uniformly
continuous in the (compact) set $[0,\tau]$. By Assumption (A.7.4),
it follows that $x_{\tau+\theta}\in
C^1([-\Delta_{\max},0];\mathbb{R}^{n})$, $\theta\in
]-\Delta_{\max},0[$, i.e. the derivative $\dot{x}_{\tau+\theta}\in
C^{0}([-\Delta_{\max},0];\mathbb{R}^{n})$. By taking into account
the Lipschitz constant $\kappa$ (computed on the bounded set
$X\times X\times U$) of function $f$, the bounds $B_{X}$ and $B_{U}$, the
following chain of inequalities holds:
\begin{equation}
\begin{array}
{lll}
 \Vert \dot{x}_{\tau+\theta} \Vert_{\infty}  & = & \sup_{\alpha\in
[-\Delta_{\max},0]} \Vert f(x(\tau+\theta+\alpha),x(\tau +\theta +\alpha-\Delta(\tau+\theta+\alpha)),\nonumber\\
& & u(\tau+\theta+\alpha-r)) \Vert \\
 & \leq & \kappa\sup_{\alpha\in [-\Delta_{\max},0]}(\Vert
x(\tau+\theta+\alpha)  \Vert  +\Vert x(\tau +\theta
+\alpha-\Delta(\tau+\theta+\alpha))\Vert  + \nonumber\\
& &
\Vert u(\tau+\theta+\alpha-r)  \Vert) \\
& \leq &  \kappa(2B_X+B_{U}), \qquad \theta\in ]-\Delta_{\max},0[.
\label{cond33}
\end{array}
\end{equation}
As far as the second derivative is concerned, the following
equality holds, for $\theta\in ]-\Delta_{\max},0[$,
\[
\begin{array}
{rcl}
\frac{d^{2}x_{\tau}(\theta)}{d\theta^{2}} & = &
J(x(\tau+\theta),x(\tau+\theta-\Delta(\tau+\theta)),u(\tau+\theta-r))\cdot\nonumber\\
& &
\left
(\begin{array}{cc}\dot{x}(\tau+\theta) \nonumber\\
\dot{x}(\tau+\theta-\Delta(\tau+\theta))(1-\dot
\Delta(\tau+\theta))
\\ 0
\end{array} \right ).\nonumber
\end{array}
\]
Since the Fr\'echet differential $J$ is continuous, since $\dot{x}_{\tau+\theta}\in C^{0}([-\Delta_{\max},0];\mathbb{R}^{n})$ and from the differentiability assumption on the state time--delay, it follows that $x_{\tau}\in PC^2([-\Delta_{\max},0];\mathbb{R}^{n})$. Moreover, by
taking into account the bound $B_{J}$ on the Fr\'echet differential $J$, the bound on the derivative $\dot{x}_{\tau+\theta}$ in (\ref{cond33}) and the definition of $M_{X}$ in (\ref{cost1}), we obtain $\Vert D^{2} x_{\tau} \Vert_{\infty}\le M_{X}$. Thus $x_{\tau}$ satisfies conditions (A.7.1) and (A.7.2). Finally by condition (A.7.3) it is readily seen that $\Vert x_{\tau}
\Vert_{\infty}\leq B_{X}^{0}$.
\end{proof}

\bigskip
The above result shows that under assumptions (A.0--1), (A.4--6) and (A.7.3--4), the regularity properties of the initial state $\xi_{0}$ in (A.7.1--2) propagate to the whole set of reachable states,
or in other words, that time--delay systems are invariant with respect to the properties in (A.7.1--2).
This is a key result that allows us to underline sufficient conditions for the existence of a countable approximation of $R_{\tau}(\Sigma)$, as formally stated hereafter.

\begin{theorem}
Consider a digital time--delay system $\Sigma$, satisfying assumptions (A.0--7). Then the set $\mathcal{A}_{X}$ defined for any $\lambda_{X}\in\mathbb{R}^{+}$ by:
\begin{equation}
\mathcal{A}_{X}(\lambda_{X})=\psi_{\lambda_{X},M_{X}}(R_{\tau}(\Sigma)),
\label{ApproxR}
\end{equation}
with $\psi_{\lambda_{X},M_{X}}$ as in (\ref{refpsi}), is a countable approximation of $R_{\tau}(\Sigma)$.
\label{coroll}
\end{theorem}

\bigskip
The proof of the above result is a direct consequence of Lemmas \ref{prop1} and \ref{prop22} and it is therefore omitted.

\subsection{Main Result}
We now have all the ingredients to define a symbolic model for digital time--delay systems. Given $\tau\in\mathbb{R}^{+}$, $\theta_{X},\theta_{D},\lambda_{U},M_{D}\in\mathbb{R}^{+}$ and $N_{X},N_{D}\in\mathbb{N}$, consider the transition system
\begin{equation}
T_{\tau,(N_{X},\theta_{X}),\lambda_{U},(N_{D},\theta_{D})}(\Sigma):=T_{\tau,\lambda_{X},\lambda_{U},\lambda_{D}}(\Sigma),
\label{sm1}
\end{equation}
where $T_{\tau,\lambda_{X},\lambda_{U},\lambda_{D}}(\Sigma)$ is defined in (\ref{sm}) with $\lambda_{X}=\Lambda(N_{X},\theta_{X},M_{X})$ and
$\lambda_{D}=\Lambda(N_{D},\theta_{D},M_{D})$.
It is readily seen that Assumptions (A.0--7) guarantees that transition system
$T_{\tau,(N_{X},\theta_{X}),\lambda_{U},(N_{D},\theta_{D})}(\Sigma)$ in (\ref{sm1}) is symbolic.
We can now present the main result of this paper.

\begin{theorem}\label{teoremaprincipale}
Consider a digital time--delay system $\Sigma
$ and any desired precision $\varepsilon\in\mathbb{R}^{+}$. Given $M_{D}\in\mathbb{R}^{+}$, suppose that assumptions (A.0--7) are satisfied. Moreover let
$\tau,\theta_{X},\theta_{D},\lambda_{U}\in\mathbb{R}^{+}$ and $N_{X},N_{D}\in\mathbb{N}$ satisfy the following inequality
\begin{equation}
\label{cond3}
\begin{array}
{l}
\max\{\beta(\varepsilon,\tau),\gamma_{U}(\lambda_{U})+ \gamma_{D}(\Lambda(N_{D},\theta_{D},M_{D}))\}+\\
\max\{\beta(\Lambda(N_{X},\theta_{X},M_{X}),\tau),\gamma_{D}(\Lambda(N_{D},\theta_{D},M_{D}))\}+
\Lambda(N_{X},\theta_{X},M_{X})\leq \varepsilon,
\end{array}
\end{equation}

with $\Lambda$ as in (\ref{lambda}) and $M_{X}$ as in (\ref{cost1}). Then transition systems $T_{\tau}(\Sigma)$ and $T_{\tau,(N_{X},\theta_{X}),\lambda_{U},(N_{D},\theta_{D})}(\Sigma)$ are A$\varepsilon$A bisimilar, i.e.
$T_{\tau,(N_{X},\theta_{X}),\lambda_{U},(N_{D},\theta_{D})}(\Sigma) \cong_{\varepsilon}^{\alt} T_{\tau}(\Sigma)$.
\label{Main}
\end{theorem}
\begin{proof}
The map $\mathcal{A}_{U}$ is a countable approximation of $U$, by Proposition \ref{DelayApprox}, the map $\mathcal{A}_{D}$ is a countable approximation of $\mathcal{D}_{\tau}$ and by Theorem \ref{coroll}, the map $\mathcal{A}_{X}$ is a countable approximation of $R_{\tau}(\Sigma)$. Choose $\lambda_{X} \in \mathbb{R}^{+}$, $\lambda_{D} \in \mathbb{R}^{+}$ and $\lambda_{U}\in\mathbb{R}^{+}$ satisfying
the inequality in (\ref{cond2}). There exist $\theta_{X}\in\mathbb{R}^{+}$ and $N_{X}\in\mathbb{N}$ so that $\lambda_{X}=\Lambda(N_{X},\theta_{X},M_{X})$, $\theta_{D}\in\mathbb{R}^{+}$ and $N_{D}\in\mathbb{N}$ so that $\lambda_{D}=\Lambda(N_{D},\theta_{D},M_{D})$, that satisfy the inequality in (\ref{cond3}). Finally the result holds as a direct application of Theorem \ref{ThMain}.
\end{proof}


\section{Construction of Symbolic Models}

The construction of the symbolic model $T_{\tau,(N_{X},\theta_{X}),\lambda_{U},(N_{D},\theta_{D})}(\Sigma)$ in (\ref{sm1}) requires the preliminary computation of countable approximations $\mathcal{A}_{U}$ and $\mathcal{A}_{D}$ of $\mathcal{U}_{\tau}$ and $\mathcal{D}_{\tau}$, respectively. While the  computation of $\mathcal{A}_{U}$ is straightforward, the computation of $\mathcal{A}_{D}$ is not so, because $\mathcal{A}_{D}$ is defined as the image through the operator $\psi_{\lambda_{D},M_{D}}$ of the set $\mathcal{D}_{\tau}$ that is composed by an infinite and uncountable number of functions. In this section we present some results that are weaker than the one in Theorem \ref{teoremaprincipale}, for which in turn, symbolic models can be effectively constructed.
The main idea is to define a suitable symbolic model $T^{*}_{\tau,(N_{X},\theta_{X}),\lambda_{U},(N_{D},\theta_{D})}(\Sigma)$ that can be effectively computed and that approximates the symbolic model $T_{\tau,(N_{X},\theta_{X}),\lambda_{U},(N_{D},\theta_{D})}(\Sigma)$, in the sense of alternating $0$--approximate simulation relation. Given
\[
T_{\tau,(N_{X},\theta_{X}),\lambda_{U},(N_{D},\theta_{D})}(\Sigma)=(Q_{2},q_{2}^{0},L_{2},\rTo_{2},O_{2},H_{2}),
\]
with $L_{2}=A_{2}\times B_{2}$, define the symbolic model:
\begin{equation}
T^{\ast}_{\tau,(N_{X},\theta_{X}),\lambda_{U},(N_{D},\theta_{D})}(\Sigma):=(Q_{2},q_{2}^{0},L^{\ast}_{2},\rTo_{2},O_{2},H_{2}),
\label{smNew}
\end{equation}
where $L^{\ast}_{2}=A_{2}\times B^{\ast}_{2}$ and
$B^{\ast}_{2}$ is the collection of all functions $b(t)=\sum_{i=0}^{N_{D}+1} \tilde{y}_{i}s_{i}(t)$,
\mbox{$t\in [-\Delta_{\max},0]$}, with $\tilde{y}_{i}\in [D]_{\theta_{D}}$. Note that by definition $B_{2}\subseteq B_{2}^{\ast}$ from which, the following result holds.
\bigskip
\begin{proposition} 
\[
T^{\ast}_{\tau,(N_{X},\theta_{X}),\lambda_{U},(N_{D},\theta_{D})}(\Sigma)\preceq_{0}^{\alt} T_{\tau,(N_{X},\theta_{X}),\lambda_{U},(N_{D},\theta_{D})}(\Sigma) \preceq_{0} T^{\ast}_{\tau,(N_{X},\theta_{X}),\lambda_{U},(N_{D},\theta_{D})}(\Sigma).
\]
\label{propT*}
\end{proposition}
The proof of the above result is a straightforward consequence of the definition of the symbolic models involved and it is therefore omitted.
We now have all the ingredients to present the following result.

\begin{theorem}\label{teoremaprincipale2}
Consider a digital time--delay system $\Sigma
$ and any desired precision $\varepsilon\in\mathbb{R}^{+}$. Given $M_{D}\in\mathbb{R}^{+}$, suppose that assumptions (A.0--7) are satisfied. Moreover let
$\tau,\theta_{X},\theta_{D},\lambda_{U}\in\mathbb{R}^{+}$ and $N_{X},N_{D}\in\mathbb{N}$ satisfy the inequality in (\ref{cond3}), with $\Lambda$ as in (\ref{lambda}) and $M_{X}$ as in (\ref{cost1}). Then:
\[
T^{\ast}_{\tau,(N_{X},\theta_{X}),\lambda_{U},(N_{D},\theta_{D})}(\Sigma)
\preceq_{\varepsilon}^{\alt}
T_{\tau}(\Sigma)
\preceq_{\varepsilon}
T^{\ast}_{\tau,(N_{X},\theta_{X}),\lambda_{U},(N_{D},\theta_{D})}(\Sigma).
\]
\end{theorem}

\begin{proof}
(Proof of $T^{\ast}_{\tau,(N_{X},\theta_{X}),\lambda_{U},(N_{D},\theta_{D})}(\Sigma) \preceq_{\varepsilon}^{\alt} T_{\tau}(\Sigma)$). By
Proposition \ref{propT*},
\begin{equation}
T^{\ast}_{\tau,(N_{X},\theta_{X}),\lambda_{U},(N_{D},\theta_{D})}(\Sigma)\preceq_{0}^{\alt} T_{\tau,(N_{X},\theta_{X}),\lambda_{U},(N_{D},\theta_{D})}(\Sigma).
\label{cond*1}
\end{equation}
By Theorem \ref{Main} and since $A\varepsilon A$ bisimulation implies $A\varepsilon A$ simulation, one gets:
\begin{equation}
T_{\tau,(N_{X},\theta_{X}),\lambda_{U},(N_{D},\theta_{D})}(\Sigma) \preceq_{\varepsilon}^{\alt} T_{\tau}(\Sigma).
\label{cond*2}
\end{equation}
Hence, by combining (\ref{cond*1}) and (\ref{cond*2}) and by a straightforward generalization of Proposition 2 in \cite{AB-TAC07}, the result follows. 
The proof of the approximate inclusion $T_{\tau}(\Sigma)\preceq_{\varepsilon} T^{\ast}_{\tau,(N_{X},\theta_{X}),\lambda_{U},(N_{D},\theta_{D})}(\Sigma)$ can be given along the lines of the proof of Theorem \ref{teoremaprincipale} and is therefore omitted. 
\end{proof}

\bigskip
From \cite{Alternating,PolaSIAM2009}, the above result guarantees that control strategies synthesized on the symbolic model $T^{\ast}_{\tau,(N_{X},\theta_{X}),\lambda_{U},(N_{D},\theta_{D})}(\Sigma)$ can be readily transferred to the original system $\Sigma$, independently of the particular realization of the time--varying delay signal $\Delta$. 
The above result is weaker than Theorem \ref{teoremaprincipale} in the sense that it does not guarantee existence of alternating approximate bisimulation between the time--delay system $\Sigma$ and the corresponding symbolic model $T^{\ast}_{\tau,(N_{X},\theta_{X}),\lambda_{U},(N_{D},\theta_{D})}(\Sigma)$ (as Theorem \ref{teoremaprincipale} does). The motivation in the introduction of the symbolic model $T^{\ast}_{\tau,(N_{X},\theta_{X}),\lambda_{U},(N_{D},\theta_{D})}(\Sigma)$ is that it can be effectively computed, as discussed hereafter.
It is easy to see that the set $B^{\ast}_{2}$ coincides with the co--domain of the operator $\psi_{\theta_{D},M_{D}}$ (see (\ref{refpsi}) and (\ref{ApproxDelay})) and it is composed by a finite number of functions. Hence, the set $B^{\ast}_{2}$ can be computed in a finite number of steps from which, the symbolic model in (\ref{smNew}) can be effectively constructed. The construction of the proposed symbolic models can be easily derived by adapting Algorithm 1 in \cite{PolaSCL10} for symbolic models of time--delay systems with constant delays to symbolic models of time--delay systems with time--varying delays. However, the adaptation of Algorithm 1 in \cite{PolaSCL10} to this framework is not efficient from the computational complexity point of view because it generally leads to large symbolic models and extensive time of computation that are not needed when solving many (symbolic) control design problems. 
A more efficient approach would construct the symbolic controller without constructing the whole symbolic model of the time--delay system. Useful insights in this direction are reported in \cite{PolaTAC12}, concerning the integrated symbolic control design of nonlinear systems. We do not report in the paper technical results generalizing the ones of \cite{PolaTAC12} to time--delay systems for lack of space. We instead illustrate in the next section, through a simple example, the computational complexity gain obtained by following this approach.

\section{An illustrative example}\label{example}
In this section we illustrate the results presented in this paper by means of a simple example. Consider the following nonlinear time--delay system:
\begin{equation}
\label{sys:example}
\Sigma:
\left\{
\begin{array}
{l}
\dot x_1(t) =-8 x_1(t)+\tanh(x_2(t-\Delta(t)))\\
\dot x_2(t) =-9 x_2(t)+\sin(x_1(t-\Delta(t)))+\cos(x_2(t)) u(t-r),
\end{array}
\right.
\end{equation}

where $t\in \mathbb{R}_0^+$, $\Delta(t)\in [\Delta_{min}, \Delta_{\max}]$ with $\Delta_{\min}= 10^{-3}$, $\Delta_{\max}=10^{-2}$, $r=10$, $d_{\min}=0.2$, and
$\tanh(x)=\frac {e^x-e^{-x}}{e^x+e^{-x}}$ for any $x\in \mathbb{R}$. We address a symbolic control design problem where the state of $\Sigma$ is requested to reach some regions of the state space within some prescribed times. More specifically, we consider the following synchronization specification:
starting from the origin, remain in the positive orthant for all times; reach the set $X_{1}=[0.01,\infty[ \times [0.1,\infty[$ in no more than $4$s, stay in the set $X_{1}$ for at least $4$s, reach the set $X_{2}=[0.08,015]\times [0.08,0.15]$ and finally remain in $X_{2}$ for at least $12$s. These requirements arise for example in multi--agent systems with shared resources in which the use of a given resource is needed to be synchronized among the agents. More complex specifications can be also considered as logics--based specification, fairness constraints, and etc. (see e.g. \cite{LTLControl}).  
In order to solve the considered control design problem we first need to check $\delta$-IDSS of $\Sigma$. 
For the system to satisfy Assumption (A.1), consider the quadratic functional (see Remark 3.9 in \cite{PepeJiangSCL}, concerning the linear increasing kernel in the second integral term) defined, for any $t\in \mathbb{R}_0^+$, $\phi =\left [\begin{array}{c}\phi_1 \\ \phi_2 \end{array}\right ]\in
C^0([-\Delta_{\max},0];\mathbb{R}^2)$, as: 
\begin{eqnarray}&&
V_{ISS}(t,\phi)=\phi^2_1(0)+\phi^2_2(0)+2\int_{-\Delta(t)}^0(\phi^2_1(\tau)+\phi^2_2(\tau))d\tau \nonumber \\ && \qquad \qquad +
\int_{-\Delta_{\max}}^0\left (
\frac{-r_{\Delta}\tau}{\Delta_{max}}+\frac{r_0(\tau+\Delta_{max})}{\Delta_{max}}\right
)\phi^2(\tau)d\tau,
\end{eqnarray}
where $r_0$, $r_{\Delta}$ are positive reals, with
$r_0>r_{\Delta}$. By applying Theorem 3.1 in \cite{PepeJiangSCL} and appropriately choosing the parameters $r_0$, $r_{\Delta}$, we obtain that $V_{ISS}$ is an ISS Lyapunov--Krasovskii functional for system (\ref{sys:example}). In particular, for $r_{\Delta}=0.2$, $r_0=0.3$, the inequality in (\ref{ineqiss}) is fulfilled with functions:
\[
\begin{array}
{lll} \beta_{ISS}(\omega,t)=2.8920 e^{-1.0870t}\omega, &
\gamma_{ISS}(\omega)= 0.9592\, \omega
&\omega,t\in\mathbb{R}^{+}_{0}.
\end{array}
\]
We now proceed with a further step and consider Assumption (A.4).
Consider the quadratic functional defined, for any $t\in
\mathbb{R}_0^+$, $\phi_1,\phi_2 \in C([-\Delta_{\max},0];\mathbb{R}^2)$, as:

\begin{eqnarray}&&
 V(t,\phi_1,\phi_2)= 
(\phi_1(0)-\phi_2(0))^T(\phi_1(0)-\phi_2(0))\nonumber\\ && \qquad +2\int_{-\Delta(t)}^{0} (\phi_1(\tau)-\phi_2(\tau))^T(\phi_1(\tau)-\phi_2(\tau))d\tau
\nonumber\\
&& \qquad +\int_{-\Delta_{\max}}^0\left(
\frac{-r_{\Delta}\tau}{\Delta_{max}}+\frac{r_0(\tau+\Delta_{max})}{\Delta_{max}}
\right )
(\phi_1(\tau)-\phi_2(\tau))^T(\phi_1(\tau)-\phi_2(\tau))d\tau, \nonumber \\ &&
\end{eqnarray}

where $r_0$, $r_{\Delta}$ are positive reals, with $r_0>r_{\Delta}$. By appropriately choosing the parameters $r_0$, $r_{\Delta}$, we obtain that $V$ is a $\delta$--IDSS Lyapunov--Krasovskii functional for system (\ref{sys:example}) and hence, by Theorem \ref{Th_KLF}, the time--delay system $\Sigma$ is $\delta$--IDSS. In particular, for $r_0=0.3$, $r_{\Delta}=0.2$, the inequality in (\ref{ineqdeltaidss}) is fulfilled with functions:

\[
\begin{array}
{llll} \beta(\omega,t)=4.3580e^{-1.0870t}\omega, &
\gamma_{U}(\omega)= 13.5647\,\omega, & \gamma_{D}(\omega)=
194.1666\,\omega, & \omega,t\in\mathbb{R}^{+}_{0}.
\end{array}
\]

Let be $B_{X}^0=0.5$, $B_{U}=0.3$, $M_{1}=0.1$, and $M_{D}=0.001$. Consequently, we obtain $B_{J}=27.9$, $\kappa=9.3$, $M_{X}=993.8845$. For a precision $\varepsilon=0.12$, we can choose $\tau=2$, $\lambda_{X}=0.02$, $\lambda_{U}=5 \cdot 10^{-4}$, and $\lambda_{D}=1.4\cdot 10^{-4}$ so that the inequality in (\ref{cond2}) is satisfied. By the definition of function $\Lambda$ in (\ref{lambda}) and since $\lambda_{X}=\Lambda(N_{X},\theta_{X},M_{X})$ and $\lambda_{D}=\Lambda(N_{D},\theta_{D},M_{D})$, one can choose $N_X=0$, $\theta_X=0.04$, $N_D=1$ and $\theta_D=6 \cdot 10^{-6}$. This choice of quantization parameters satisfies the conditions in (A.7.3) and (\ref{cond3}). By generalizing the algorithms proposed in \cite{PolaSCL10}, the estimated time to construct the symbolic model $T^{\ast}_{\tau,(N_{X},\theta_{X}),\lambda_{U},(N_{D},\theta_{D})}(\Sigma)$ is about $203,215s$. 
Since the expected time of computation is rather high, in the sequel we adapt the algorithms proposed in \cite{PolaTAC12} towards the efficient symbolic control design of nonlinear time--delay systems. The symbolic control strategy obtained is reported hereafter, where $(n_{1},n_{2})\rTo^{u}(n_{1}^{+},n_{2}^{+})$ stands for $((n_{1}\theta_{X},n_{1}\theta_{X}),(n_{2}\theta_{X},n_{2}\theta_{X}))\rTo^{u\theta_{U}}((n_{1}^{+}\theta_{X},n_{1}^{+}\theta_{X}),$ $(n_{2}^{+}\theta_{X},n_{2}^{+}\theta_{X}))$:
\begin{equation}
\label{control}
\begin{array}
{ccccccc}
(0,0) & \rTo^{186}  & (4,30) & \rTo^{-396}  & (3,22) & \rTo^{248}  & \\

(4,31) & \rTo^{-562}  & (3,20) & \rTo^{-268}  & (3,24) & \rTo^{-546}  & \\

(2,20) &  \rTo^{-484}  & (3,21) & \rTo^{388}  & (4,33) & \rTo^{234}  & \\

(4,31) & \rTo^{-220}  & (3,25) & \rTo^{542}  & (4,35) & \rTo^{-560}  & \\

(3,19) & \rTo^{-74}  & (3,27) & \rTo^{-142}  & (3,26). \\& & & & & &
\end{array}
\end{equation}

\begin{figure}[ht]
 \centering
 {\includegraphics[scale=0.5]{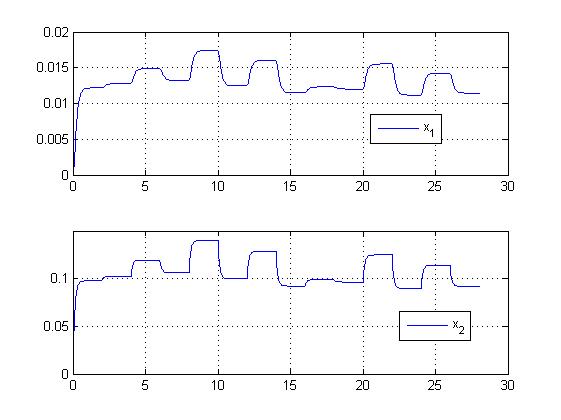}}
 \caption{Evolution of the state variables of $\Sigma$ with initial condition $\xi_{0}=0$, time--varying delay signal in (\ref{TDsignal}), and the  control strategy reported in (\ref{control}).}
 \label{figs}
 \end{figure}

The running time needed for solving the given symbolic control design problem is $7,692s$ by using a laptop with CPU Intel Core 2 Duo T5500 @ $1.66$GHz. 
Figure \ref{figs} shows the evolution of the state variables of $\Sigma$ with the unknown time--delay signal:
\begin{equation}
\label{TDsignal}
\Delta(t)=\frac{\Delta_{\max}+\Delta_{\min}}{2} + \frac{\Delta_{\max}-\Delta_{\min}}{2} \sin(0.01\,t),\,\,t\in\mathbb{R}^{+}_{0}.
\end{equation}
It is readily seen that the synchronization specification is indeed satisfied.

%

\section{Conclusion}
In this paper we studied existence of symbolic models for nonlinear control systems with time--varying and unknown time--delay signals. We first introduced the notion of $\delta$--IDSS 
and provided a characterization of this property, in terms of Lyapunov--Krasovskii functionals. We then showed that $\delta$--IDSS time--delay systems admit symbolic models that are alternating approximately bisimilar to the original systems with a precision that can be rendered as small as desired. We finally presented a symbolic model that is shown to be an alternating approximate simulation of the original system and that can be computed in a finite number of steps under a boundedness assumption on the sets of states, inputs and delays of the system. \\

\bigskip
\textbf{Acknowledgement:} The authors thank Prof. Paulo Tabuada (University of California at Los Angeles, USA) for fruitful discussions on the topic of this paper and Dr. Alessandro Borri (University of L'Aquila, Italy) for the software implementation of the algorithms used in Section 7.

\bibliographystyle{alpha}
\bibliography{biblio1}

\bigskip
\section{Appendix: Notation} \label{sec:notation}

The symbols $\mathbb{N}$, $\mathbb{Z}$, $\mathbb{R}$, $\mathbb{R}^{+}$ and $\mathbb{R}_{0}^{+}$ denote the sets of natural, integer, real, positive and nonnegative real numbers, respectively. Given a vector $x\in\mathbb{R}^{n}$ the $i$--th element of $x$ is denoted by $x_{i}$; furthermore $\Vert x\Vert$
denotes the infinity norm of $x$; we recall that \mbox{$\Vert x\Vert:=\max\{|x_1|,|x_2|,...,|x_n|\}$}, where $|x_i|$ is the absolute value of $x_i$. For a given $x\in \mathbb {R}^n$ and a given $s\in \mathbb{R}^{+}$, the symbol $\mathcal {B}_{s}(x)$ denotes the set $\{y\in \mathbb {R}^n: \Vert y-x\Vert \le s\}$.
For any $A\subseteq \mathbb{R}^{n}$ and \mbox{$\theta\in{\mathbb{R}^{+}}$} define
\begin{equation}
[A]_{\theta}:=\{a\in A\,\,|a_{i}=k_{i}\theta,\,\,\,k_{i}\in\mathbb{Z},i=1,...,n\}.
\label{grid}
\end{equation}
Given a measurable and locally essentially bounded function \mbox{$f:\mathbb{R}_{0}^{+}\rightarrow\mathbb{R}^{n}$}, the \mbox{(essential)} supremum norm of $f$ is denoted by $\Vert f\Vert_{\infty}$; we recall that \mbox{$\Vert f\Vert_{\infty}:=(ess)\sup\{\Vert f(t)\Vert,$ $t\geq0\}$}. For a given time $\tau\in\mathbb{R}^{+}$, define $f_{\tau}$ so that $f_{\tau}(t)=f(t)$, for any $t\in [0,\tau[$, and $f(t)=0$ elsewhere; $f$ is said to be locally essentially bounded if for any $\tau\in\mathbb{R}^{+}$, $f_{\tau}$ is essentially bounded.
A continuous function $\gamma:\mathbb{R}_{0}^{+}%
\rightarrow\mathbb{R}_{0}^{+}$ is said to belong to class $\mathcal{K}$ if it is strictly increasing and \mbox{$\gamma(0)=0$}; $\gamma$ is said to belong to class $\mathcal{K}_{\infty}$ if \mbox{$\gamma\in\mathcal{K}$} and $\gamma(r)\rightarrow \infty$ as $r\rightarrow\infty$. A continuous function \mbox{$\beta:\mathbb{R}_{0}^{+}\times\mathbb{R}_{0}^{+}\rightarrow\mathbb{R}_{0}^{+}$} is said to belong to class $\mathcal{KL}$ if for each fixed $s$, the map $\beta(r,s)$ belongs to class $\mathcal{K}$ with respect to $r$ and, for each fixed $r$, the map $\beta(r,s)$ is decreasing with respect to $s$ and $\beta(r,s)\rightarrow0$ as \mbox{$s\rightarrow\infty$}. Given $k,n\in\mathbb{N}$ with $n\geq 1$ and $I=[a,b]\subseteq \mathbb{R}$, $a,b\in\mathbb{R}$, $a<b$ let $C^{k}(I;\mathbb{R}^{n})$ be the space of functions $f:I\rightarrow \mathbb{R}^{n}$ that are continuously differentiable $k$ times. Given $k\geq 1$, let $PC^{k}(I;\mathbb{R}^{n})$ be the space of $C^{k-1}(I;\mathbb{R}^{n})$ functions $f:I\rightarrow \mathbb{R}^{n}$ whose $k$--th derivative exists except in a finite number of reals, and it is bounded, i.e. there exist $\gamma_{0},\gamma_{1},...,\gamma_{s}\in\mathbb{R}^{+}$ with  $a=\gamma_{0}<\gamma_{1}<...<\gamma_{s}=b$ so that $D^{k}\,f$ is defined on each open interval $(\gamma_{i},\gamma_{i+1})$, $i=0,1,...,s-1$ and $\max_{i=0,1,...,s-1}\sup_{t\in(\gamma_{i},\gamma_{i+1})}\Vert D^{k}\,f(t) \Vert _{\infty}<\infty$. For any continuous function $x(s)$, defined on $-\Delta\leq s<a$, $a>0$, and any fixed $t$, $0\leq t<a$, the standard symbol $x_{t}$ will denote the element of $C^{0}([-\Delta,0];\mathbb{R}^{n})$ defined by $x_{t}(\theta)=x(t+\theta)$, $-\Delta\leq\theta\leq 0$.
The identity map on a set $A$ is denoted by $1_{A}$. Given two sets $A$ and $B$, if $A$ is a subset of $B$ we denote by \mbox{$\imath_{A}:A\hookrightarrow B$} or simply by $\imath$ the natural inclusion map taking any $a\in A$ to \mbox{$\imath (a) = a \in B$}. Given a function $f:A\rightarrow B$ the symbol $f(A)$ denotes the image of $A$ through $f$, i.e. $f(A):=\{b\in B:\exists a\in A$ s.t.
$b=f(a)\}$. If $C\subseteq A$, the symbol $f|_{C}$ denotes the restriction of function $f$ to $C$, i.e. $f|_{C}(c)=f(c)$ for any $c\in C$. Given a relation $R\subseteq A\times B$, $R^{-1}$ denotes the inverse relation of $R$, i.e. $R^{-1}:=\{(b,a)\in B\times A:( a,b)\in A\times B\}$.

\end{document}